\documentclass[%
  a4paper,
  onecolumn,
  colorlinks,
]{preprint}

\usepackage[english]{babel}

\usepackage{graphicx}
\usepackage{subcaption}
\usepackage{pgfplotstable}
\usepackage{booktabs}

\usepackage{tikz}
\usepackage{pgfplots}
  \pgfplotsset{compat = 1.18}
  \usepgfplotslibrary{external}
  \tikzset{external/system call = {%
    pdflatex \tikzexternalcheckshellescape
      -halt-on-error
      -interaction=batchmode
      -jobname "\image" "\texsource"}}
  \tikzexternaldisable

\usepackage{amsmath}
\usepackage{amssymb}
\usepackage{amsthm}

\usepackage[linesnumbered, lined, algoruled]{algorithm2e}

\usepackage{todonotes}

\theoremstyle{plain}\newtheorem{theorem}{Theorem}
\theoremstyle{plain}
\theoremstyle{definition}\newtheorem{remark}{Remark}

\newcommand{\R}{\ensuremath{\mathbb{R}}}
\newcommand{\C}{\ensuremath{\mathbb{C}}}

\newcommand{\bA}{\ensuremath{\boldsymbol{A}}}
\newcommand{\bB}{\ensuremath{\boldsymbol{B}}}
\newcommand{\bC}{\ensuremath{\boldsymbol{C}}}
\newcommand{\bD}{\ensuremath{\boldsymbol{D}}}
\newcommand{\bE}{\ensuremath{\boldsymbol{E}}}

\newcommand{\bI}{\ensuremath{\boldsymbol{I}}}

\newcommand{\bK}{\ensuremath{\boldsymbol{K}}}
\newcommand{\bL}{\ensuremath{\boldsymbol{L}}}
\newcommand{\bM}{\ensuremath{\boldsymbol{M}}}
\newcommand{\bN}{\ensuremath{\boldsymbol{N}}}

\newcommand{\bQ}{\ensuremath{\boldsymbol{Q}}}
\newcommand{\bR}{\ensuremath{\boldsymbol{R}}}
\newcommand{\bS}{\ensuremath{\boldsymbol{S}}}
\newcommand{\bT}{\ensuremath{\boldsymbol{T}}}
\newcommand{\bU}{\ensuremath{\boldsymbol{U}}}
\newcommand{\bV}{\ensuremath{\boldsymbol{V}}}
\newcommand{\bW}{\ensuremath{\boldsymbol{W}}}
\newcommand{\bX}{\ensuremath{\boldsymbol{X}}}

\newcommand{\bc}{\ensuremath{\boldsymbol{c}}}
\newcommand{\bff}{\ensuremath{\boldsymbol{f}}}
\newcommand{\bt}{\ensuremath{\boldsymbol{t}}}
\newcommand{\bu}{\ensuremath{\boldsymbol{u}}}
\newcommand{\bv}{\ensuremath{\boldsymbol{v}}}

\newcommand{\bx}{\ensuremath{\boldsymbol{x}}}
\newcommand{\by}{\ensuremath{\boldsymbol{y}}}

\newcommand{\bPsi}{\ensuremath{\boldsymbol{\Psi}}}

\newcommand{\bLcal}{\ensuremath{\boldsymbol{\mathcal{L}}}}
\newcommand{\bPcal}{\ensuremath{\boldsymbol{\mathcal{P}}}}

\newcommand{\herm}{\ensuremath{\mathsf{H}}}
\newcommand{\trans}{\ensuremath{\mkern-1.5mu\mathsf{T}}}
\DeclareMathOperator{\real}{Re}
\DeclareMathOperator{\imag}{Im}
\DeclareMathOperator{\mspan}{span}
\DeclareMathOperator{\mdiag}{diag}
\DeclareMathOperator*{\argmin}{argmin}
\DeclareMathOperator*{\argmax}{argmax}

\newcommand{\eye}[1]{\ensuremath{\bI_{#1}}}
\newcommand{\iunit}{\ensuremath{\mathfrak{i}}}
\newcommand{\fro}{\operatorname{F}}
\newcommand{\Hinf}{\ensuremath{\mathcal{H}_{\infty}}}
\newcommand{\matlab}{\mbox{MATLAB}}

\newcommand{\plotfontsize}{\small}

\definecolor{matlabblue}{HTML}{0072BD}
\definecolor{matlaborange}{HTML}{D95319}
\definecolor{matlabyellow}{HTML}{EDB120}
\definecolor{matlabpurple}{HTML}{7E2F8E}
\definecolor{matlabgreen}{HTML}{77AC30}
\definecolor{matlablightblue}{HTML}{4DBEEE}
\definecolor{matlabred}{HTML}{A2142F}

\newcommand{\adi}{\texttt{ADI}}
\newcommand{\tadi}{\texttt{TADI}}
\newcommand{\aditrunc}{\texttt{ADI (trunc)}}
\newcommand{\taditrunc}{\texttt{TADI (trunc)}}
\newcommand{\tadirand}{\texttt{TADI(rand)}}
\newcommand{\tadieig}{\texttt{TADI(eig)}}

\tikzstyle{adi_imp_style} = [
  matlabblue,
  solid,
  line width = 2pt
]

\tikzstyle{adi_exp_style} = [
  matlabpurple,
  dashed,
  line width = 2.5pt
]

\tikzstyle{eigtadi_imp_style} = [
  matlaborange,
  solid,
  line width   = 2pt,
  mark options = {solid},
  mark         = *
]

\tikzstyle{eigtadi_exp_style} = [
  matlabgreen,
  dashed,
  line width   = 2.5pt,
  mark options = {solid},
  mark         = *
]

\tikzstyle{tadi_rand_style} = [
  matlabgreen,
  solid,
  line width   = 2pt,
  mark options = {solid},
  mark         = triangle*
]

\tikzstyle{eigtadi_fullheur} = [
  matlabblue,
  dotted,
  line width   = 2pt,
  mark options = {solid},
  mark         = square*
]

\tikzstyle{eigtadi_resheur} = [
  matlabgreen,
  dashed,
  line width   = 2pt,
  mark options = {solid},
  mark         = diamond*
]

\tikzstyle{adi_trunc_imp_style} = [
  matlabgreen,
  solid,
  line width = 2pt,
  mark       = triangle*
]

\tikzstyle{eigtadi_trunc_imp_style} = [
  matlabpurple,
  solid,
  line width   = 2pt,
  mark options = {solid},
  mark         = diamond*
]

\begin{document}

\title{A tangential low-rank ADI method for solving indefinite
  Lyapunov equations}

\author[$\ast$]{Rudi Smith}
\affil[$\ast$]{Department of Mathematics,
  Virginia Tech, Blacksburg, VA 24061, USA.\authorcr
  \email{smithrgh@vt.edu}, \orcid{0000-0002-6592-2482}}

\author[$\dagger$]{Steffen W. R. Werner}
\affil[$\dagger$]{Department of Mathematics,
  Division of Computational Modeling and Data Analytics, and
  National Security Institute,
  Virginia Tech, Blacksburg, VA 24061, USA.\authorcr
  \email{steffen.werner@vt.edu}, \orcid{0000-0003-1667-4862}
}

\shorttitle{Tangential ADI for indefinite Lyapunov equations}
\shortauthor{R. Smith, S. W. R. Werner}
\shortdate{2026-09-13}
\shortinstitute{}


\abstract{
Continuous-time algebraic Lyapunov equations have become an essential tool in
  various applications.
  In the case of large-scale sparse coefficient matrices and indefinite
  constant terms, indefinite low-rank factorizations have successfully been used
  to allow methods like the alternating direction implicit (ADI) iteration to
  efficiently compute accurate approximations to the solution of the Lyapunov
  equation.
  However, classical block-type approaches quickly increase in computational
  costs when the rank of the constant term grows.
  While the regular truncation of solution approximations during the
  iteration may be a remedy, it does not resolve the issue of computing
  high-dimensional updates and may even compromise the numerical accuracy
  of the approximation.
  In this paper, we propose a novel tangential reformulation of the
  ADI iteration that preserves the exactness of the iterative updates while 
  allowing for the efficient construction of low-rank approximations to the
  solution of Lyapunov equations with indefinite right-hand sides even in the
  case of constant terms with higher ranks.
  We provide adaptive methods for the selection of the corresponding ADI
  parameters, namely shifts and tangential directions, which allow for the
  automatic application of the method to any relevant problem setting.
  The effectiveness of the developed algorithms is illustrated in several
  numerical examples.
}

\novelty{We propose a novel tangential ADI method for the solution of Lyapunov
equations with indefinite constant terms.
For the selection of shift parameters and tangential directions, we develop
adaptive algorithms, which update the ADI parameters during the iteration
as needed.}

\maketitle


\section{Introduction}%
\label{sec:intro}

Continuous-time algebraic Lyapunov equations are one of the most
important foundational tools in systems and control theory with various
applications in model order reduction~\cite{BenMS05, BenS17, Gre88, Moo81} and
controller design~\cite{Arm75, HeiKW24, Son98}.
They also often appear as the backbone of differential and nonlinear matrix
equation solvers~\cite{BenS13, LanMS15, MenOPetal18, SaaW24}.
In general, continuous-time algebraic Lyapunov equations can be written in
the form
\begin{equation} \label{eqn:lyap}
  \bA \bX \bE^{\herm} + \bE \bX \bA^{\herm} + \bB \bR \bB^{\herm} = 0,
\end{equation}
with the coefficient matrices $\bA, \bE \in \C^{n \times n}$,
$\bB \in \C^{n \times m}$, and
$\bR \in \C^{m \times m}$ Hermitian.
Here, we denote by $\bM^{\herm} = \overline{\bM}^{\trans}$ the
transposed conjugate of the matrix $\bM \in \C^{n \times m}$.
The task is to find a Hermitian matrix $\bX \in \C^{n \times n}$, which solves
the Lyapunov equation~\cref{eqn:lyap}.
In this work, we consider the classical assumption that the matrix pencil
$\lambda \bE - \bA$ is Hurwitz, meaning that all its finite eigenvalues lie in
the open left half-plane.
Also, for simplicity of illustration, we assume in the following the case that
$\bE$ is invertible;
however, we outline modifications for the case of non-invertible
$\bE$ matrices later in \Cref{sec:singularE}.

A variety of methods is available for the solution of~\cref{eqn:lyap}.
In the case of dense small-scale coefficient matrices, direct solvers as well
as high-performing iterative solvers have been
developed~\cite{BarS72, BenCQ98, Ham82, Rob80}.
In the case when $n$ becomes large ($\gtrapprox 10^{4}$),
the matrices $\bA$ and $\bE$ typically appear to be only sparsely populated.
However, the solution $\bX$ remains a dense matrix, which complicates
computations in terms of time and memory constraints.
With the additional assumption that the number of columns of $\bB$ is small
and $\bR$ is positive definite, new solvers have been developed, which
do not interfere with the sparse structure of $\bA$ and $\bE$ while
approximating the solution $\bX$ via low-rank
factors~\cite{LiW02, BenLP08, HocS95, DruS11, Sim07}.
Recently, an extension of the Alternating Direction Implicit (ADI) method has
been developed for the case that $\bR$ is allowed to be
indefinite~\cite{LanMS14, LanMS15}.

One of the main assumptions in solvers for large-scale sparse Lyapunov
equations is that the constant term has a low rank $m$.
This not only promotes the fast decay of the singular values of the solution
$\bX$, which allows its approximation via low-rank factors in the first place,
the factors $\bB$ and $\bR$ are actively involved in the construction of the
approximation~\cite{BenLP08, DruS11, LanMS14, Sim07}.
However, it has been observed that there are
applications~\cite{HeiKW24, ShaSS15}, in
which the constant term is not low-rank and cannot be well approximated in
low-rank form, but the solution $\bX$ still exhibits a fast singular value decay
for which it can be well approximated by low-rank factors.
A remedy for this problem is the framework of tangential Krylov methods, where
underlying subspaces and solution approximations constructed in
the methods are compressed along suitable directions.
For the special case of $\bR = \eye{m}$ in~\cref{eqn:lyap},
tangential rational Krylov methods have been discussed in~\cite{DruSZ14}, while
a corresponding ADI approach has been proposed in~\cite{WolPL16}.
We note that the assumption $\bR = \eye{m}$ is essential for the results
presented in~\cite{DruSZ14, WolPL16}.
In particular, the computational formulas, the convergence theory and the
proposed selection methods for tangential directions do not hold
in the presence of an arbitrary Hermitian center matrix $\bR$
in~\cref{eqn:lyap} as we will show later in \Cref{sec:tangdirect}.

In this work, we propose a new tangential ADI method for the solution of
large-scale sparse Lyapunov equations with indefinite constant terms.
This method allows the efficient approximation of the solution to Lyapunov
equations via low-rank factors even in the case of high-rank constant terms.
In contrast to truncation-based approaches that compress periodically solely the
solution approximation, for example, via a truncated SVD, the proposed method
is theoretically exact and preserves the relation between the solution
approximation, the residual terms, and the computed update factors.
Furthermore, the proposed method resolves the issue of computing low-rank
update terms, which cannot be handled via truncation otherwise.
Due to its practical importance, we provide theoretical results and
algorithms separately for the case of Lyapunov equations with real coefficients
of the form
\begin{equation} \label{eqn:lyap_real}
  \bA \bX \bE^{\trans} + \bE \bX \bA^{\trans} + \bB \bR \bB^{\trans} = 0,
\end{equation}
with $\bA, \bE \in \R^{n \times n}$, $\bB \in \R^{n \times m}$, and
$\bR \in \R^{m \times m}$ symmetric.
In contrast to~\cref{eqn:lyap}, computational steps have to be modified so
that the constructed low-rank factors preserve the realness of the coefficients.
To answer the fundamental questions about choosing ADI shifts and the tangential
directions, we develop two adaptive algorithms that compute shifts and
directions during the ADI iteration as needed.

The rest of the manuscript is organized as follows:
In \Cref{sec:prelim}, we summarize the classical indefinite
factorized ADI method for the general case of Lyapunov equations with complex
coefficients. Based on these formulas and computational steps, we then develop
the tangential ADI method for Lyapunov equations with complex and real
coefficients in \Cref{sec:tangADI}.
The adaptive generation of suitable ADI shift parameters and the selection
of tangential directions are discussed in \Cref{sec:param}.
In \Cref{sec:numerics}, the proposed methods are tested in numerical
experiments verifying their effectiveness in comparison to the
current state-of-the-art approach.
The paper is concluded in \Cref{sec:conclusions}.


\section{Indefinite factorized ADI method}%
\label{sec:prelim}

The tangential approach, which we propose in this work, is based on the
low-rank ADI iteration for Lyapunov equations of the form~\cref{eqn:lyap}. Since the extension of low-rank ADI methods to indefinite right-hand
sides has been previously established in~\cite{LanMS14, LanMS15}, we briefly
summarize here the principal idea and the resulting indefinite factorized ADI
scheme.
We formulate the algorithm for complex and real coefficient matrices as these
will serve as the standard approaches for numerical
comparisons later in \Cref{sec:numerics}.


\subsection{Basic iteration scheme}%
\label{sec:basicADI}

Originally, the Alternating Direction Implicit (ADI) method was
proposed for the numerical solution of parabolic and elliptic differential
equations~\cite{PeaR55}, known to converge superlinearly under mild
assumptions~\cite{BecG10}.
The extension of the ADI method to Lyapunov equations of the
form~\cref{eqn:lyap} with real coefficient matrices was first proposed
in~\cite{LanMS14}.

Following the detailed derivations in~\cite{LanMS15}
and~\cite{Kue16}, we adapt here the iteration for complex matrices.
For the solution of~\cref{eqn:lyap}, it is well documented that a low-rank
constant term encourages a fast singular value decay of the solution 
matrix $\bX$ so that it can be well approximated via low-rank matrix
factors~\cite{BakES15, Pen00, TruV07}.
Applying the indefinite symmetric structure from the constant term
in~\cref{eqn:lyap} to the ADI iterates so that
$\bX_{j} = \bL_{j} \bD_{j} \bL_{j}^{\herm}$ and assuming a zero
initial value $\bX_{0} = 0$ results in the basic low-rank factorized ADI
iteration given by
\begin{subequations} \label{eqn:LjDj_final}
\begin{align} \label{eqn:LjDj_final_Vj}
    \bV_{j} & = (\bA + \alpha_{j} \bE)^{-1} (\bA - \overline{\alpha_{j-1}} \bE)
      \bV_{j - 1}, \\
    \bL_{j} & = \begin{bmatrix} \bL_{j - 1} & \bV_{j} \end{bmatrix},\\
    \bD_{j} & = -2 \mdiag( \real(\alpha_{1}), \ldots, \real(\alpha_{j}) )
      \otimes \bR,
\end{align}
\end{subequations}
for $j = 2, 3, \ldots$, and with the initial matrices
$\bV_{1} = (\bA + \alpha_{1} \bE)^{-1} \bB$ and
$\bL_{1} = \bV_{1}$.
We note that in implementations, dimensionless empty matrices
$\bL_{0} = [~]$ and $\bR_{0} = [~]$ are used for the initial iterate since
$\bL_{0}$ and $\bR_{0}$ do not contribute to the solution approximation. It follows from~\cref{eqn:LjDj_final} that the ADI iterates
$\bX_{j}$ are constructed as the sum of rank-$m$ matrices
$\bX_{j} = \sum_{k = 1}^{j} -2\real(\alpha_{k}) \bV_{k} \bR \bV_{k}^{\herm}$.


\subsection{Implicit residual formulation}%
\label{sec:impResADI}

The residual of ADI iterates $\bX_{j} = \bL_{j} \bD_{j} \bL_{j}^{\herm}$ is
generally given by
\begin{equation} \label{eqn:expRes}
  \bLcal_{j} = \bA \bX_{j} \bE^{\herm} + \bE \bX_{j} \bA^{\herm} +
    \bB \bR \bB^{\herm},
\end{equation}
which ideally converges during the iteration to zero.
To effectively compute the residual during the iteration as well as to improve
the computational costs of the ADI steps in~\cref{eqn:LjDj_final}, we will
make use of the fact that in the ADI iteration, the residual~\cref{eqn:expRes}
can be implicitly assembled in the same factorized form as the constant term
in~\cref{eqn:lyap} and our solution iterates
$\bX_{j} = \bL_{j} \bD_{j} \bL_{j}^{\herm}$.

As established in~\cite{BenKS13b}, the ADI Lyapunov residual can be
reformulated via recursive updates using the ADI iterates.
Adapted for the indefinite factorized ADI method~\cref{eqn:LjDj_final} with
the initial zero iterate $\bX_{0} = 0$, the residual at any step $j \geq 1$ has
at most rank $m$ and is given by $\bLcal_{j} = \bW_{j} \bR \bW_{j}^{\herm}$.
Thereby, the matrix $\bW_{j} \in \C^{n \times m}$ follows the recursion formula
\begin{equation} \label{eqn:resW}
  \bW_{j} = (\bA - \overline{\alpha_{j}} \bE) \bV_{j}
    = \bW_{j - 1} - 2\real(\alpha_{j}) \bE \bV_{j},
\end{equation}
with the initial residual factor $\bW_{0} = \bB$.

The factorized structure of the Lyapunov residual in~\cref{eqn:resW}
allows for the efficient computation of residual norms during the ADI
iteration.
This can be used to set appropriate stopping criteria of the general form
$\lVert \bLcal_{j} \rVert < \epsilon \lVert \bB \bR \bB^{\herm} \rVert$,
where $\epsilon > 0$ is a user-defined tolerance and $\lVert . \rVert$
denotes a suitable matrix norm.
The resulting ADI iteration for Lyapunov equations~\cref{eqn:lyap}
is shown in \Cref{alg:LDLADIcplx}.

\begin{algorithm}[t]
  \SetAlgoHangIndent{1pt}
  \DontPrintSemicolon
  
  \caption{Indefinite factorized ADI method: complex case.}
  \label{alg:LDLADIcplx}
  
  \KwIn{Coefficient matrices $\bA, \bE, \bB, \bR$ from~\cref{eqn:lyap},
    shift parameters $\alpha_{1}, \ldots, \alpha_{\ell} \in \C_{-}$,
    residual tolerance $\epsilon > 0$.}
  \KwOut{Approximate solution factors $\bL_{j - 1}, \bD_{j - 1}$ so that
    $\bL_{j - 1} \bD_{j - 1} \bL_{j - 1}^{\herm} \approx \bX$.}
  
  Initialize $\bW_{0} = \bB$, $\bL_{0} = [~]$, and $j = 1$.\;
  
  \While{$\lVert \bW_{j - 1} \bR \bW_{j - 1}^{\herm} \rVert \geq
    \epsilon \lVert \bB \bR \bB^{\herm} \rVert$
    \label{alg:LDLADIcplx_loop}}{
    Solve $(\bA + \alpha_{j} \bE) \bV_{j} = \bW_{j - 1}$ for $\bV_{j}$.\;
    
    Update the residual factor
      $\bW_{j} = \bW_{j - 1} - 2\real(\alpha_{j}) \bE \bV_{j}$.\;
    
    Update the solution factor $\bL_{j} =
      \begin{bmatrix} \bL_{j - 1} & \bV_{j} \end{bmatrix}$.\;
    
    Increment $j \gets j + 1$.\;
  }

  Compute the solution factor
    $\bD_{j - 1} = -2 \mdiag(\real(\alpha_{1}), \ldots, \real(\alpha_{j - 1}))
    \otimes \bR$.\;
\end{algorithm}

In the case that $m < n$, the cost of computing the residual can be
further reduced in the spectral or Frobenius norms.
Both of these norms allow for computationally efficient reformulations that are
based on the fact that the non-zero spectra of
$\bLcal_{j} = \bW_{j} \bR \bW_{j}^{\herm} \in \C^{n \times n}$
and $\bW_{j}^{\herm} \bW_{j} \bR \in \C^{m \times m}$ are identical.
For the spectral norm, we have that
\begin{equation} \label{eqn:twonorm_tmp}
  \lVert \bLcal_{j} \rVert_{2} =
    \lVert \bW_{j} \bR \bW_{j}^{\herm} \rVert_{2} =
      \max\limits_{\lambda \in \Lambda(\bW_{j}^{\herm} \bW_{j} \bR)}
      \lvert \lambda \rvert,
\end{equation}
while for the Frobenius norm, it holds that
\begin{equation} \label{eqn:fronorm_tmp}
  \lVert \bLcal_{j} \rVert_{\fro} =
    \lVert \bW_{j} \bR \bW_{j}^{\herm} \rVert_{\fro} =
    \sqrt{\sum\limits_{k = 1}^{m} \lvert \lambda_{k} \rvert^{2}},
    \quad\text{where}\quad
    \{ \lambda_{1}, \ldots, \lambda_{m} \} =
      \Lambda(\bW_{j}^{\herm} \bW_{j} \bR).
\end{equation}
While this reordering works well in the case that $\bR = \eye{m}$,
we note that due to the general loss of symmetry of the matrix
$\bW_{j}^{\herm} \bW_{j} \bR$ when $\bR \neq \eye{m}$, numerical issues
may occur.
A suitable alternative to~\cref{eqn:twonorm_tmp} and~\cref{eqn:fronorm_tmp}
can be given using the economy-size QR factorization of the outer residual
matrices, $\bW_{j} = \bQ_{j} \bU_{j}$, so that the residual can also be written
as $\bLcal_{j} = \bQ_{j} (\bU_{j} \bR \bU_{j}^{\herm}) \bQ_{j}^{\herm}$.
Thereby, the inner matrix given by $\bU_{j} \bR \bU_{j}^{\herm}$ is symmetric
with dimensions $m$.
Due to the unitary invariance of the spectral and Frobenius norms, the residual
norms can be computed via
\begin{equation} \label{eqn:two_fro_norms}
  \lVert \bLcal_{j} \rVert_{2} = \lVert \bU_{j} \bR \bU_{j}^{\herm} \rVert_{2}
  \quad\text{and}\quad
  \lVert \bLcal_{j} \rVert_{\fro} = \lVert \bU_{j} \bR \bU_{j}^{\herm}
  \rVert_{\fro}.
\end{equation}
We refer the reader to~\cite[Lem.~3.1]{DenKLetal25} for the proof of unitary
norm invariance with rectangular unitary matrices and further discussions about
these norm computations.
These two expressions in~\cref{eqn:two_fro_norms} can be used to efficiently
replace the norm computations done in \Cref{alg:LDLADIcplx_loop} of \Cref{alg:LDLADIcplx}.


\subsection{Real coefficients and real solution factors}%
\label{sec:realADI}

In many important use cases, the  coefficient matrices in the Lyapunov
equation are real~\cref{eqn:lyap_real}.
Consequently, the solution $\bX$ becomes real and symmetric, and in the
construction of low-rank solutions, we wish the solution factors to preserve
this realness, too.
In the case that the shift parameters $\alpha_{1}, \ldots, \alpha_{\ell}$
in \Cref{alg:LDLADIcplx} are chosen to be real, the algorithm naturally
constructs real solution factors $\bL_{j}$ and $\bD_{j}$ where $\bD_{j}$
is symmetric.
However, for most coefficient matrices, complex shifts are needed to ensure
fast convergence of the ADI method~\cite{Wac13, BakES15}.

To guarantee realness when complex shifts are utilized, it is
standard practice to combine two consecutive ADI steps with complex conjugate
shifts, $\alpha_{j+1} = \overline{\alpha_{j}}$.
As shown in the literature for standard formulations, this conjugate double ADI
update allows the update solution and residual factors to be computed in real
arithmetic, and halving the number of linear systems that must be solved.
Applied to the indefinite factorization, tracking the real and imaginary
components of $\bV_{j}$ yields the basis change formulation
\begin{subequations} \label{eqn:basischange}
\begin{align}
  \bX_{j + 1} & = \bL_{j - 1} \bD_{j - 1} \bL_{j - 1}^{\trans} +
    \begin{bmatrix} \bV_{j} & \bV_{j+1} \end{bmatrix}
    \begin{bmatrix} -2 \real(\alpha_{j}) \bR & 0 \\ 0 & -2 \real(\alpha_{j})
    \bR \end{bmatrix}
    \begin{bmatrix} \bV_{j}^{\trans} \\ \bV_{j + 1}^{\trans} \end{bmatrix} \\
  & = \bL_{j - 1} \bD_{j - 1} \bL_{j - 1}^{\trans} +
    \begin{bmatrix} \real(\bV_{j}) & \imag(\bV_{j}) \end{bmatrix}
    \bT_{j} \bPsi_{j} \bT_{j}^{\herm}
    \begin{bmatrix} \real(\bV_{j})^{\trans} \\ \imag(\bV_{j})^{\trans}
    \end{bmatrix},
\end{align}
\end{subequations}
with $\delta_{j} = \real(\alpha_{j}) / \imag(\alpha_{j})$ and the matrices
\begin{equation*}
  \bT_{j} = \begin{bmatrix} \eye{m} & \eye{m} \\ \iunit \eye{m} &
    (2 \delta_{j} - \iunit) \eye{m} \end{bmatrix}
  \quad\text{and}\quad
  \bPsi_{j} = \begin{bmatrix} -2 \real(\alpha_{j}) \bR & 0 \\ 0 &
    -2 \real(\alpha_{j}) \bR \end{bmatrix}.
\end{equation*}
Alternatively, the center term can be expressed via real matrices as
\begin{equation*}
   \bT_{j} \bPsi_{j} \bT_{j}^{\herm} = \bM_{j} \bPsi_{j} \bM_{j}^{\trans},
   \quad\text{with}\quad
   \bM_{j} = \sqrt{2} \begin{bmatrix} \eye{m} & 0 \\ \delta_{j} \eye{m} &
     \sqrt{\delta_{j}^{2} + 1} \eye{m} \end{bmatrix}.
\end{equation*}
Integrating the matrix $\bM_{j}$ into the update for the outer solution factors
leads to the real-valued ADI update scheme
\begin{align*}
  \bL_{j + 1} & = \begin{bmatrix} \bL_{j - 1} &
    \begin{bmatrix} \real(\bV_{j}) & \imag(\bV_{j}) \end{bmatrix} \bM_{j}
    \end{bmatrix} \\
  & = \begin{bmatrix} \bL_{j - 1} & \sqrt{2} \big( \real(\bV_{j}) +
    \delta_{j} \imag(\bV_{j}) \big) & \sqrt{2 (\delta_{j}^{2} + 1)}
    \imag(\bV_{j}) \end{bmatrix}, \\
  \bD_{j + 1} & = \begin{bmatrix} \bD_{j - 1} & & \\
    & -2 \real(\alpha_{j}) \bR & \\ & & -2 \real(\alpha_{j}) \bR \end{bmatrix}.
\end{align*}
The resulting algorithm for Lyapunov equations with real coefficient
matrices~\cref{eqn:lyap_real} is outlined in \Cref{alg:LDLADIreal}.

\begin{algorithm}[t]
  \SetAlgoHangIndent{1pt}
  \DontPrintSemicolon
  
  \caption{Indefinite factorized ADI method: real case.}
  \label{alg:LDLADIreal}
  
  \KwIn{Real coefficient matrices $\bA, \bE, \bB, \bR$
    from~\cref{eqn:lyap_real},
    shift parameters $\alpha_{1}, \ldots, \alpha_{\ell} \in \C_{-}$ closed
    under conjugation and ordered as pairs,
    residual tolerance $\epsilon > 0$.}
  \KwOut{Approximate solution factors $\bL_{j - 1}, \bD_{j - 1}$ so that
    $\bL_{j - 1} \bD_{j - 1} \bL_{j - 1}^{\trans} \approx \bX$.}
  
  Initialize $\bW_{0} = \bB$, $\bL_{0} = [~]$, $\eta = \sqrt{2}$, and $j = 1$.\;
  
  \While{$\lVert \bW_{j - 1} \bR \bW_{j - 1}^{\trans} \rVert \geq
    \epsilon \lVert \bB \bR \bB^{\trans} \rVert$}{
    
    Solve $(\bA + \alpha_{j} \bE) \bV_{j} = \bW_{j - 1}$ for $\bV_{j}$.\;
    
    \eIf{$\alpha_{j}$ is real}{
      Update the residual and the solution factors via
        \vspace{-.5\baselineskip}
        \begin{align*}
          \bW_{j} & = \bW_{j - 1} - 2 \alpha_{j} \bE \bV_{j}, \\
          \bL_{j} & = \begin{bmatrix} \bL_{j - 1} & \bV_{j} \end{bmatrix}.
        \end{align*}\;
        \vspace{-1.5\baselineskip}
    }{
      Compute the scaling factor
        $\delta_{j} = \real(\alpha_{j}) / \imag(\alpha_{j})$.\;

      Update the residual and the solution factors via
        \vspace{-.5\baselineskip}
        \begin{align*}
          \bW_{j + 1} & = \bW_{j - 1} - 4 \real(\alpha_{j}) \bE \big(
            \real(\bV_{j}) + \delta_{j} \imag(\bV_{j}) \big), \\
          \bL_{j + 1} & = \begin{bmatrix} \bL_{j - 1} &
            \eta \big( \real(\bV_{j}) + \delta_{j} \imag(\bV_{j}) \big) &
            \eta \sqrt{\delta_{j}^{2} + 1} \imag(\bV_{j}) \end{bmatrix}.
        \end{align*}\;
        \vspace{-1.5\baselineskip}

      Increment $j \gets j + 1$.\;
    }

  Increment $j \gets j + 1$.\;
  }

  Compute the solution factor
    $\bD_{j - 1} = -2 \mdiag(\real(\alpha_{1}), \ldots, \real(\alpha_{j - 1}))
    \otimes \bR$.\;
\end{algorithm}


\section{Tangential indefinite low-rank ADI method}%
\label{sec:tangADI}

In this section, we introduce the concept of compressing solution factor updates
in the ADI method from \Cref{sec:prelim} using tangential directions.
We develop the theoretical foundations for a residual formulation of the
method, and we provide the algorithms for the cases of complex as well as 
real coefficient matrices.
The special case of non-invertible $\bE$ matrices is discussed at the end.


\subsection{Tangential updates and the residual formula}%
\label{sec:tangTheory}

While generally an effective tool to solve large-scale Lyapunov equations,
the classical ADI method as described in \Cref{sec:prelim} operates blockwise.
In other words, every ADI step extends the size of the solution factor
$\bL_{j}$ by $m$ columns and the block diagonal term $\bD_{j}$ by $m$ rows
and columns, where $m$ is the rank of the right-hand side term
in~\cref{eqn:lyap} and~\cref{eqn:lyap_real}.
Thus, the ADI steps become computationally demanding in terms of memory for
storing the solution factors as well as computational costs for solving
linear systems with $m$ right-hand sides.
A remedy to this, originating in matrix interpolation theory~\cite{BalGR90},
is the projection of updates $\bV_{j} \in \C^{n \times m}$ onto suitable
tangential directions $\bt_{j} \in \C^{m}$ so that the new updates take the
form $\bv_{j} = \bV_{j} \bt_{j} \in \C^{n}$.
Thereby, only linear systems with a single right-hand side have to be solved in
every ADI step, and the solution factors grow only by a single column and
diagonal entry, respectively.

For the special case of~\cref{eqn:lyap_real} with $\bR = \eye{m}$, a tangential
version of the ADI method has been proposed previously in~\cite{WolPL16}.
The derivation therein is based on the observation that the update terms in
the ADI method build a basis for a rational Krylov subspace.
Since tangential interpolation is necessary in the construction of optimal
interpolations of rational functions via rational Krylov
subspaces~\cite{GalVV04}, the idea in~\cite{WolPL16} was to potentially
transfer those results towards the ADI method.
In contrast, our work aims for the general reduction of computational costs in
the ADI method using the concept of tangential compression, diverging
from the rational Krylov framework in~\cite{WolPL16} in several key aspects.
First, regarding the problem setting, we consider the general case of complex
coefficient matrices in~\cref{eqn:lyap} with a potentially indefinite constant
terms, whereas~\cite{WolPL16} is restricted to real coefficients and a positive
semi-definite right-hand side of the form $\bB\bB^{\trans}$.
Second, the indefinite center matrix $\bR$ in~\cref{eqn:lyap} drastically
alters the convergence requirements:
while the method in~\cite{WolPL16} can utilize arbitrary normalized tangential
directions, we empirically show that indefinite problems strictly require
directions derived from the eigenvectors of $\bR$ to prevent divergence of the
ADI.
Finally, to ensure practical efficiency and robust convergence, we introduce an
adaptive tangential direction selection that complements the popular projection
shift strategy, overcoming the parameter selection limitations present in
previous tangential approaches.

The following theorem considers the compression of the ADI update terms
$\bV_{j}$ from~\cref{eqn:LjDj_final} down to rank $1$ using tangential
directions and provides a formula for the corresponding Lyapunov
residual~\cref{eqn:expRes}.

\begin{theorem}[Tangential ADI update]%
  \label{thm:tangADI}
  Consider the Lyapunov equation~\cref{eqn:lyap} and assume that $\bR$ is
  invertible.
  Let $\bX_{j - 1} = \bL_{j - 1} \bD_{j - 1} \bL_{j - 1}^{\herm}$ be the ADI
  iterate at step $j-1$, with the associated Lyapunov residual
  $\bLcal_{j - 1} = \bW_{j - 1} \bR \bW_{j - 1}^{\herm}$.
  Given the next ADI shift $\alpha_{j} \in \C_{-}$ and the tangential direction
  $\bt_{j} \in \C^{m}$, the rank-$1$ compressed ADI update of the form
  \begin{equation*}
    \bX_{j} = \bX_{j - 1} + d_{j} \bv_{j}  \bv_{j}^{\herm}
  \end{equation*}
  is given via the scalar
  \begin{equation*}
    d_{j} = \frac{-2 \real(\alpha_{j})}{\bt_{j}^{\herm} \bR^{-1} \bt_{j}},
  \end{equation*}
  and $\bv_{j} \in \C^{n}$ as the solution to the linear system
  \begin{equation*}
    (\bA + \alpha_{j} \bE) \bv_{j} = \bW_{j - 1} \bt_{j}.
  \end{equation*}
  The corresponding Lyapunov residual has the form
  $\bLcal_{j} = \bW_{j} \bR \bW_{j}^{\herm}$, where
  \begin{equation*}
    \bW_{j} = \bW_{j - 1} - \frac{2 \real(\alpha_{j}) \bE \bv_{j}
      \bt_{j}^{\herm} \bR^{-1}}{\bt_{j}^{\herm} \bR^{-1} \bt_{j}}.
  \end{equation*}
\end{theorem}
\begin{proof}
  To prove the results, we make use of the connection between the ADI iterates
  and the implicit residual expression in~\cref{eqn:resW}.
  We proceed by equating the expression for the residual given
  above with the Lyapunov residual in~\cref{eqn:expRes}.
  Substituting the solution update
  $\bX_{j} = \bX_{j - 1} + \bv_{j} d_{j} \bv_{j}^{\herm}$
  into~\cref{eqn:expRes} yields
  \begin{align*}
    \bLcal_{j} & = (\bA \bX_{j - 1} \bE^{\herm} + \bE \bX_{j - 1} \bA^{\herm} +
      \bB \bR \bB^{\herm}) + d_{j} (\bA \bv_{j} \bv_{j}^{\herm} \bE^{\herm} +
      \bE \bv_{j} \bv_{j}^{\herm} \bA^{\herm}) \\
    & = \bLcal_{j - 1} + d_{j} \left( (\bA \bv_{j}) (\bE \bv_{j})^{\herm} +
      (\bE \bv_{j}) (\bA \bv_{j})^{\herm} \right).
  \end{align*}
  From the definition of $\bv_{j}$, it holds that
  $\bA \bv_{j} = \bW_{j - 1} \bt_{j} - \alpha_{j} \bE \bv_{j}$.
  Inserting this into the previous expression for the residual and using
  the identity $\alpha_{j} + \overline{\alpha_{j}} = 2 \real(\alpha_{j})$,
  we have that
  \begin{subequations}
  \begin{align}
    \bLcal_{j} & = \bLcal_{j - 1} + d_{j}
      \left( (\bW_{j - 1} \bt_{j} - \alpha_{j} \bE \bv_{j})
      (\bE \bv_{j})^{\herm} + (\bE \bv_{j}) (\bW_{j - 1} \bt_{j} -
      \alpha_{j} \bE \bv_{j})^{\herm} \right) \\
    \label{eqn:tangADI_tmp1}
    & = \bLcal_{j - 1} + d_{j} \left( (\bW_{j - 1} \bt_{j})
      (\bE \bv_{j})^{\herm} + (\bE \bv_{j}) (\bW_{j - 1} \bt_{j})^{\herm} -
      2 \real(\alpha_{j}) (\bE \bv_{j}) (\bE \bv_{j})^{\herm} \right).
  \end{align}
  \end{subequations}
  On the other side, consider the residual given by
  $\bLcal_{j} = \bW_{j} \bR \bW_{j}^{\herm}$, and let
  the residual factor follow the factorized structure
  in~\cref{eqn:resW} with a general rank-$1$ update so that
  $\bW_{j} = \bW_{j - 1} + \bE \bv_{j} \bc_{j}^{\herm}$ holds, with an
  unknown vector $\bc_{j} \in \C^{m}$.
  Inserting the rank-$1$ update into the factorized residual then yields
  \begin{subequations}
  \begin{align}
    \bLcal_{j} & = (\bW_{j - 1} + \bE \bv_{j} \bc_{j}^\herm) \bR
      (\bW_{j - 1} + \bE \bv_{j} \bc_{j}^{\herm})^{\herm} \\
    \label{eqn:tangADI_tmp2}
    & = \bLcal_{j - 1} + (\bW_{j - 1} \bR \bc_{j}) (\bE \bv_{j})^{\herm} +
      (\bE \bv_{j}) (\bW_{j - 1} \bR \bc_{j})^{\herm} +
      (\bc_{j}^{\herm} \bR \bc_{j}) (\bE \bv_{j}) (\bE \bv_{j})^{\herm}. 
  \end{align}
  \end{subequations}
  To match up the residuals, we compare the coefficients of the terms in
  \cref{eqn:tangADI_tmp1} and~\cref{eqn:tangADI_tmp2}.
  Thereby, we find that
  \begin{equation*}
    d_{j} \bW_{j - 1} \bt_{j} = \bW_{j - 1} \bR \bc_{j}.
  \end{equation*}
  Under the assumption that $\bW_{j - 1}$ has full column rank, multiplying
  both sides of this equation by the pseudoinverse of $\bW_{j - 1}$ leads to
  $\bc_{j} = d_{j} \bR^{-1} \bt_{j}$.
  We note that this $\bc_{j}$ is always a particular solution to the equation,
  also in the cases when $\bW_{j - 1}$ has not full column rank.
  By comparing the scalar coefficients of the
  $(\bE \bv_{j}) (\bE \bv_{j})^{\herm}$ term, we find that
  \begin{equation*}
    -2 \real(\alpha_{j}) d_{j} = \bc_{j}^{\herm} \bR \bc_{j}.
  \end{equation*}
  With the expression for $\bc_{j}$ from above, it follows
  \begin{equation*}
    -2 \real(\alpha_{j}) d_{j} =
      (d_{j} \bR^{-1} \bt_{j})^{\herm} \bR (d_{j} \bR^{-1} \bt_{j}) =
      d_{j}^{2} (\bt_{j}^{\herm} \bR^{-1} \bR \bR^{-1} \bt_{j}) =
      d_{j}^{2} (\bt_{j}^{\herm} \bR^{-1} \bt_{j}).
  \end{equation*}
  Under the assumption that $d_{j} \neq 0$, this equation can be resolved
  to
  \begin{equation*}
    d_{j} = \frac{-2 \real(\alpha_{j})}{\bt_{j}^{\herm} \bR^{-1} \bt_{j}}.
  \end{equation*}
  Finally, substituting this expression for $d_{j}$ back into the rank-$1$
  update for the residual gives the desired factor update
  \begin{align*}
    \bW_{j} & = \bW_{j - 1} + \bE \bv_{j} \bc_{j}^{\herm} \\
    & = \bW_{j - 1} + d_{j} \bE \bv_{j} (\bR^{-1} \bt_{j})^{\herm} \\
    & = \bW_{j - 1} - \frac{2 \real(\alpha_{j}) \bE \bv_{j}
      \bt_{j}^{\herm} \bR^{-1}}{\bt_{j}^{\herm} \bR^{-1} \bt_{j}},
  \end{align*}
  which concludes the proof.
\end{proof}

We note that in contrast to the classical implicit residual
formulation as in~\cref{eqn:resW}, the tangential ADI given by \Cref{thm:tangADI} has the additional assumption that
the center matrix $\bR$ in~\cref{eqn:lyap} is invertible.
This can always be satisfied by truncating any rank-deficient parts from
the constant term in~\cref{eqn:lyap}.
Such a ``rank truncation'' is also strongly recommended as a
pre-processing step for the classical block ADI method
(\Cref{sec:prelim}) to eliminate unnecessary computational costs
resulting from non-contributing parts of the constant term.

Following the results of \Cref{thm:tangADI}, we can derive the
approximate solution factors $\bL_{j}$ and $\bD_{j}$ in a similar fashion as
done in the block case in \Cref{sec:prelim}.
Let the previous ADI iterate be given in factorized form
$\bX_{j - 1} = \bL_{j - 1} \bD_{j - 1} \bL_{j - 1}{\herm}$, then
the next iterate that includes the rank-$1$ update from \Cref{thm:tangADI} is given by
\begin{equation*}
  \bX_{j} = \bX_{j - 1} + d_{j} \bv_{j} \bv_{j}^{\herm} =
    \begin{bmatrix} \bL_{j - 1} & \bv_{j} \end{bmatrix}
    \begin{bmatrix} \bD_{j - 1} & 0 \\ 0 & d_{j} \end{bmatrix}
    \begin{bmatrix} \bL_{j - 1}^{\herm} \\ \bv_{j}^{\herm} \end{bmatrix}.
\end{equation*}
The resulting indefinite factorized tangential ADI method based on the updating
scheme from \Cref{thm:tangADI} is summarized in \Cref{alg:LDLTADIcplx}.

\begin{algorithm}[t]
  \SetAlgoHangIndent{1pt}
  \DontPrintSemicolon
  
  \caption{Indefinite factorized tangential ADI method: complex case.}
  \label{alg:LDLTADIcplx}
  
  \KwIn{Coefficient matrices $\bA, \bE, \bB, \bR$ from~\cref{eqn:lyap},
    shift parameters $\alpha_{1}, \ldots, \alpha_{\ell} \in \C_{-}$,
    tangential directions $\bt_{1}, \ldots, \bt_{\ell} \in \C^{m}$,
    residual tolerance $\epsilon > 0$.}
  \KwOut{Approximate solution factors $\bL_{j-1}, \bD_{j-1}$ so that
    $\bL_{j-1} \bD_{j-1} \bL_{j-1}^{\herm} \approx \bX$.}
  
  Initialize $\bW_{0} = \bB$, $\bL_{0} = [~]$, and $j = 1$.\;
  
  \While{$\lVert \bW_{j - 1} \bR \bW_{j - 1}^{\herm} \rVert \geq
    \epsilon \lVert \bB \bR \bB^{\herm} \rVert$}{
    
    Solve $(\bA + \alpha_{j} \bE) \bv_{j} = \bW_{j - 1} \bt_{j}$
      for $\bv_{j}$.\;

    Update the residual factor
      $\bW_{j} = \bW_{j - 1} - \frac{2 \real(\alpha_{j})\bE \bv_{j}
        \bt_{j}^{\herm} \bR^{-1}}{\bt_{j}^{\herm} \bR^{-1} \bt_{j}}$.\;
    
    Update the solution factor
      $\bL_{j} = \begin{bmatrix} \bL_{j - 1} & \bv_{j} \end{bmatrix}$.
    
    Increment $j \gets j + 1$.\;
  }

  Compute the solution factor
      $\bD_{j - 1} = \mdiag \left(
        \frac{-2 \real(\alpha_{1})}{\bt_{1}^{\herm} \bR^{-1} \bt_{1}}, \ldots,
        \frac{-2 \real(\alpha_{j - 1})}{\bt_{j - 1}^{\herm} \bR^{-1}
        \bt_{j - 1}} \right)$.\;
\end{algorithm}

We note that \Cref{alg:LDLTADIcplx} corresponds to the block approach
in \Cref{alg:LDLADIcplx}.
The main differences between the two algorithms come from the restriction of
the ADI updates to rank $1$, which leads to the solution of linear systems with
only a single right-hand side and the inclusion of an update vector into the
solution factor in \Cref{alg:LDLTADIcplx} rather than a matrix with
$m$-columns as in \Cref{alg:LDLADIcplx}.
Since the residual has the same structure in both algorithms, we may use the
same techniques for the computationally efficient evaluation as described in \Cref{sec:impResADI}.

Finally, the similarity of the proposed tangential method to the classical
approach allows us to make some general conclusions about the convergence
of \Cref{alg:LDLTADIcplx}.
By replacing the tangential directions in \Cref{alg:LDLTADIcplx} by
identity matrices $\eye{m}$, we directly recover the classical approach.
Therefore, for a suitable choice of tangential directions
$\bt_{1}, \ldots, \bt_{\ell}$ and shifts $\alpha_{1}, \ldots, \alpha_{\ell}$,
the tangential ADI can recover the solution factors from the block ADI in \Cref{alg:LDLADIcplx}.
However, we note that not all choices of tangential directions will
yield convergence for the ADI method.
This will be further discussed in \Cref{sec:tangdirect}.


\subsection{Real-valued tangential ADI updates}%
\label{sec:tangReal}

As in the classical block ADI, we typically wish to construct real solution
factors $\bL_{j}$ and $\bD_{j}$ when the Lyapunov equation has real coefficient
matrices~\cref{eqn:lyap_real}.
Similar to the results in \Cref{sec:realADI}, we can make use of the
structure resulting in consecutive ADI steps when using complex conjugate
shifts.
The resulting real-valued update scheme is shown in the following theorem.

\begin{theorem}[Tangential conjugate double ADI update]%
  \label{thm:tangConjUpdate}
  Consider the Lyapunov equation with real coefficient matrices
  in~\cref{eqn:lyap_real} and assume $\bR$ is invertible.
  Let $\bW_{j - 1}$ be a real-valued residual factor, let the
  next two ADI shifts be complex conjugates so that $\alpha_{j}$ and
  $\alpha_{j + 1} = \overline{\alpha_{j}}$, and let the corresponding
  tangential directions be identical and real so that $\bt_{j} \in \R^{m}$
  and $\bt_{j + 1} = \bt_{j}$.
  Then, the iterate after two steps is real-valued,
  $\bX_{j + 1} \in \R^{n \times n}$, and can be written via a rank-$2$ update
  of the form
  \begin{equation*}
    \bX_{j + 1} = \bX_{j - 1} +
      \begin{bmatrix} \bff_{j} & \bff_{j+1} \end{bmatrix}
      \begin{bmatrix} d_{j} & 0 \\ 0 & d_{j} \end{bmatrix}
      \begin{bmatrix} \bff_{j}^{\trans} \\ \bff_{j+1}^{\trans} \end{bmatrix},
  \end{equation*}
  where the update components are given by
  \begin{align*}
    \bff_{j} & = \sqrt{2} \big( \real(\bv_{j}) + \delta_{j}
      \imag(\bv_{j}) \big), \\
    \bff_{j + 1} & = \sqrt{2} \sqrt{\delta_{j}^{2} + 1} \imag(\bv_{j}),\\
    d_{j} & = \frac{-2\real(\alpha_{j})}{\bt_{j}^{\trans} \bR^{-1} \bt_{j}},
  \end{align*}
  with the scaling factor $\delta_{j} = \real(\alpha_{j}) / \imag(\alpha_{j})$
  and the vector
  $\bv_{j} = (\bA + \alpha_{j} \bE)^{-1} \bW_{j - 1}\bt_{j}$.
  The corresponding residual factor is given by
  \begin{equation*}
    \bW_{j + 1} = \bW_{j - 1} - \frac{4 \real(\alpha_{j}) \bE
      \big( \real(\bv_{j}) + \delta_{j} \imag(\bv_{j}) \big) \bt_{j}^{\trans}
      \bR^{-1}}{\bt_{j}^{\trans} \bR^{-1} \bt_{j}}.
  \end{equation*}
\end{theorem}
\begin{proof}
  Following \Cref{thm:tangADI}, the update to the solution approximation
  after two consecutive tangential ADI steps is in general given via
  \begin{equation*}
    \bX_{j+1} = \bL_{j - 1} \bD_{j - 1} \bL_{j - 1}^{\trans} +
      \begin{bmatrix} \bv_{j} & \bv_{j + 1} \end{bmatrix}
      \begin{bmatrix} \frac{-2 \real(\alpha_{j})}{\bt_{j}^{\trans}
      \bR^{-1} \bt_{j}} & 0 \\ 0 &
      \frac{-2 \real(\alpha_{j})}{\bt_{j}^{\trans} \bR^{-1} \bt_{j}}
      \end{bmatrix} \begin{bmatrix} \bv_{j}^{\herm}\\\bv_{j + 1}^{\herm}
      \end{bmatrix},
  \end{equation*}
  where the update vectors are the solutions to the two linear systems
  \begin{equation*}
    (\bA + \alpha_{j} \bE) \bv_{j} = \bW_{j - 1} \bt_{j}
    \quad\text{and}\quad
    (\bA + \overline{\alpha_{j}} \bE) \bv_{j + 1} = \bW_{j} \bt_{j}.
  \end{equation*}
  Since all the coefficient matrices as well as the previous residual factor
  $\bW_{j - 1}$ and the tangential direction $\bt_{j}$ are real, we can apply
  the results from the standard conjugate double step scheme and the
  basis change in~\cref{eqn:basischange}, which allow us to rewrite the
  solution update by splitting the update vector $\bv_{j}$ into its real and
  imaginary parts so that
  \begin{equation*}
    \bX_{j + 1} = \bL_{j - 1} \bD_{j - 1} \bL_{j - 1}^{\trans} +
      \begin{bmatrix} \real(\bv_{j}) & \imag(\bv_{j}) \end{bmatrix}
      \bM_{j} \bPsi_{j} \bM_{j}^{\trans}
      \begin{bmatrix} \real(\bv_{j})^{\trans} \\ \imag(\bv_{j})^{\trans}
      \end{bmatrix},
  \end{equation*}
  where the center matrices are given by
  \begin{equation*}
    \bM_{j} = \sqrt{2} \begin{bmatrix} 1 & 0 \\ \delta_{j} & 
      \sqrt{\delta_{j}^2 + 1} \end{bmatrix}
    \quad\text{and}\quad
    \bPsi_{j} = \begin{bmatrix} \frac{-2 \real(\alpha_{j})}{\bt_{j}^{\trans}
      \bR^{-1} \bt_{j}} & 0 \\ 0 & \frac{-2 \real(\alpha_{j})}%
      {\bt_{j}^{\trans} \bR^{-1} \bt_{j}} \end{bmatrix},
  \end{equation*}
  with the scaling factor $\delta_{j} = \real(\alpha_{j}) / \imag(\alpha_{j})$.
  Integrating the transformation $\bM_{j}$ into the outer factors via
  \begin{align*}
    \begin{bmatrix} \real(\bv_{j}) & \imag(\bv_{j}) \end{bmatrix} \bM_{j} & =
      \begin{bmatrix} \sqrt{2} \left(\real(\bv_{j}) + \delta_{j}
      \imag(\bv_{j}) \right) & \sqrt{2} \sqrt{\delta_{j}^{2} + 1}
      \imag(\bv_{j}) \end{bmatrix}\\
    & = \begin{bmatrix} \bff_{j} & \bff_{j + 1} \end{bmatrix}
  \end{align*}
  and labeling the diagonal elements of $\Psi_{j}$ as
  $d_{j} = \frac{-2 \real(\alpha_{j})}{\bt_{j}^{\trans} \bR^{-1} \bt_{j}}$,
  yields the desired result for the real-valued update of the
  solution approximation.

  For the update of the residual factor, we have from \Cref{thm:tangADI}
  that
  \begin{equation*}
    \bW_{j} = \bW_{j - 1} - \frac{2 \real(\alpha_{j}) \bE \bv_{j}
      \bt_{j}^{\trans} \bR^{-1}}{\bt_{j}^{\trans} \bR^{-1} \bt_{j}}
  \end{equation*}
  holds, with the update vector $\bv_{j}$ from above.
  By observing that all coefficient matrices as well as the previous
  residual factor $\bW_{j - 1}$ and the tangential direction $\bt_{j}$
  are real, we may split the update vector into its real and imaginary
  part $\bv_{j} = \real(\bv_{j}) + \iunit \imag(\bv_{j})$, which satisfy
  the systems of linear equations
  \begin{align*}
    \big( \bA + \real(\alpha_{j}) \bE \big) \real(\bv_{j}) -
      \imag(\alpha_{j}) \bE \imag(\bv_{j}) & = \bW_{j - 1} \bt_{j}, \\
    \big( \bA + \real(\alpha_{j}) \bE \big) \imag(\bv_{j}) +
      \imag(\alpha_{j}) \bE \real(\bv_{j}) & = 0.
  \end{align*}
  Rearranging the second equation in terms of $\real(\bv_{j})$ allows us to
  rewrite the residual update as \begin{equation*}
    \bW_{j} = \bW_{j - 1} + \frac{2 \real(\alpha_{j})
      \Big( \big(\bA + \real(\alpha_{j}) \bE \big) \imag(\bv_{j}) -
      \iunit \imag(\alpha_{j}) \bE \imag(\bv_{j}) \Big)
      \bt_{j}^{\trans} \bR^{-1}}{\imag(\alpha_{j})\bt_{j}^{\trans} \bR^{-1} \bt_{j}}.
  \end{equation*}
  This allows us to write the update vector alternatively via
  \begin{align*}
    \bv_{j + 1} & = (\bA + \overline{\alpha_{j}} \bE)^{-1} \bW_{j} \bt_{j} \\
    & = (\bA + \overline{\alpha_{j}} \bE)^{-1} \bW_{j - 1} \bt_{j} +
      2 \delta_{j} (\bA + \overline{\alpha_{j}} \bE)^{-1}\\
    & \quad{}\times{}
      \Big( \big( \bA + \real(\alpha_{j}) \bE \big) \imag(\bv_{j}) -\iunit
      \imag(\alpha_{j}) \bE \imag(\bv_{j}) \Big)\\
    &= \overline{\bv_{j}} + 2 \delta_{j} \imag(\bv_{j}).
  \end{align*}
  Finally, for the corresponding residual factor we get
  \begin{align*}
    \bW_{j + 1} & = \bW_{j} - \frac{2 \real(\alpha_{j}) \bE \bv_{j + 1}
      \bt_{j}^{\trans} \bR^{-1}}{\bt_{j}^{\trans} \bR^{-1} \bt_{j}} \\
    & = \bW_{j - 1}  - \frac{2 \real(\alpha_{j}) \bE \bv_{j}
      \bt_{j}^{\trans} \bR^{-1} + 2 \real(\alpha_{j}) \bE \big(
      \overline{\bv_{j}} + 2 \delta_{j} \imag(\bv_{j}) \big)
      \bt_{j}^{\trans} \bR^{-1}}{\bt_{j}^{\trans} \bR^{-1} \bt_{j}} \\
    & = \bW_{j - 1}  - \frac{4 \real(\alpha_{j}) \bE \big( \real(\bv_{j}) +
      \delta_{j} \imag(\bv_{j}) \big) \bt_{j}^{\trans} \bR^{-1}}%
      {\bt_{j}^{\trans} \bR^{-1} \bt_{j}},
  \end{align*}
  which concludes the proof.
\end{proof}

Two important points in \Cref{thm:tangConjUpdate} are the assumptions
that the tangential direction $\bt_{j}$ is a real vector and that the same
direction is used for both complex conjugate shifts.
These assumptions provide the necessary structure to make \Cref{thm:tangConjUpdate} the tangential equivalent of the
standard real-valued block ADI double step.
If one were to choose complex conjugate directions
$\bt_{j+1} = \overline{\bt_{j}}$, the intermediate residual factor
$\bW_{j}$ would become complex.
Consequently, the algebraic recovery of the second update vector in the proof
of \Cref{thm:tangConjUpdate} via
$\bv_{j + 1} = (\bA + \overline{\alpha_{j}} \bE)^{-1} \bW_{j} \bt_{j+1}$ fails
to simplify into the combination of the real terms $\real(\bv_{j})$ and
$\imag(\bv_{j})$.
This would force the solution of an additional full-dimensional linear system,
negating the computational benefits of the double step.
Additionally, we do not see these assumptions to be restrictive in any capacity
considering the role of the tangential directions.
Since in the case of real coefficient matrices, these are used to essentially
compress only real-valued terms in the ADI method, there is no reason to assume
that a complex direction may yield more information about the solution.
Furthermore, we will see later in \Cref{sec:param} that these assumptions will
naturally be satisfied in our proposed method for selecting tangential
directions.

The tangential ADI method for Lyapunov equations with real coefficients
based on \Cref{thm:tangConjUpdate} is outlined in \Cref{alg:LDLTADIreal}.
As before, we can observe strong similarities between the tangential approach
in \Cref{alg:LDLTADIreal} and its classical block equivalent in \Cref{alg:LDLADIreal}.
The differences as desired are the restriction of update terms to
scalars and vectors rather than $m$-dimensional matrices.
Also, we would like to remind the reader that the residual follows the
same factorized structure as before so that the techniques mentioned in \Cref{sec:impResADI} should be used for efficient residual norm computations.

\begin{algorithm}[t]
  \SetAlgoHangIndent{1pt}
  \DontPrintSemicolon
  
  \caption{Indefinite factorized tangential ADI method: real case.}
  \label{alg:LDLTADIreal}

  \KwIn{Real coefficient matrices $\bA, \bE, \bB, \bR$
    from~\cref{eqn:lyap_real},
    shift parameters $\alpha_{1}, \ldots, \alpha_{\ell} \in \C_{-}$ closed
    under conjugation and ordered as pairs,
    real tangential directions $\bt_{1}, \ldots, \bt_{\ell} \in \R^{m}$,
    residual tolerance $\epsilon > 0$.}
  \KwOut{Approximate solution factors $\bL_{j - 1}, \bD_{j - 1}$ so that
    $\bL_{j - 1} \bD_{j - 1} \bL_{j - 1}^{\trans} \approx \bX$.}
  
  Initialize $\bW_{0} = \bB$, $\bL_{0} = [~]$, $\eta = \sqrt{2}$, and $j = 1$.\;
  
  \While{$\lVert \bW_{j - 1} \bR \bW_{j - 1}^{\herm} \rVert \geq
    \epsilon \lVert \bB \bR \bB^{\herm} \rVert$}{
    
    Solve $(\bA + \alpha_{j} \bE) \bv_{j} = \bW_{j - 1} \bt_{j}$
      for $\bv_{j}$.\;

    \eIf{$\alpha_{j}$ is real}{
      Update the residual and the solution factors via
        \vspace{-.5\baselineskip}
        \begin{align*}
          \bW_{j} & = \bW_{j - 1} - \frac{2 \alpha_{j} \bE \bv_{j}
            \bt_{j}^{\trans} \bR^{-1}}{\bt_{j}^{\trans} \bR^{-1} \bt_{j}}, \\
          \bL_{j} & = \begin{bmatrix} \bL_{j - 1} & \bv_{j} \end{bmatrix}.
        \end{align*}\;
        \vspace{-1.5\baselineskip}
    }{
      Compute the scaling factor
        $\delta_{j} = \real(\alpha_{j}) / \imag(\alpha_{j})$.\;

      Update the residual and the solution factors via
        \vspace{-.5\baselineskip}
        \begin{align*}
          \bW_{j + 1} & = \bW_{j - 1} - \frac{4 \real(\alpha_{j}) \bE \big(
            \real(\bv_{j}) + \delta_{j} \imag(\bv_{j}) \big)
            \bt_{j}^{\trans} \bR^{-1}}{\bt_{j}^{\trans} \bR^{-1} \bt_{j}}, \\
          \bL_{j + 1} & = \begin{bmatrix} \bL_{j - 1} &
            \eta \big( \real(\bv_{j}) + \delta_{j} \imag(\bv_{j}) \big) &
            \eta \sqrt{\delta_{j}^{2} + 1} \imag(\bv_{j}) \end{bmatrix}.
        \end{align*}\;
        \vspace{-1.5\baselineskip}

      Increment $j \gets j + 1$.\;
    }

    Increment $j \gets j + 1$.\;
  }

  Compute the solution factor
    $\bD_{j - 1} = \mdiag \left(
    \frac{-2 \real(\alpha_{1})}{\bt_{1}^{\trans} \bR^{-1} \bt_{1}}, \ldots,
    \frac{-2 \real(\alpha_{j - 1})}{\bt_{j - 1}^{\trans} \bR^{-1}
    \bt_{j - 1}} \right)$.\;
\end{algorithm}


\subsection{\texorpdfstring{%
  Non-invertible $\bE$ matrices and projected Lyapunov equations}%
  {Non-invertible E matrices and projected Lyapunov equations}}%
\label{sec:singularE}

For the ease of presentation, we assumed so far that the $\bE$ matrix
in~\cref{eqn:lyap} and~\cref{eqn:lyap_real} is invertible.
However, in certain applications, for example, those involving
dynamical systems~\cite{BenS17}, the $\bE$ matrix becomes
non-invertible.
Under the assumption that the matrix pencil $\lambda \bE - \bA$ is regular,
i.e., there exists a $\lambda \in \C$ so that $\det(\lambda E - A) \neq 0$,
one typically considers the solution of~\cref{eqn:lyap}
and~\cref{eqn:lyap_real} over the subspace of finite eigenvalues of
$\lambda \bE - \bA$.
This subspace-restricted problem is given via the projected Lyapunov
equation
\begin{subequations} \label{eqn:projLyap}
\begin{align}
  \bA \bX \bE^{\herm} + \bE \bX \bA^{\herm} +
    \bPcal_{\operatorname{\ell}} \bB \bR \bB^{\herm}
    \bPcal_{\operatorname{\ell}}^{\herm} & = 0, \\
  \bPcal_{\operatorname{r}} \bX \bPcal_{\operatorname{r}}^{\herm} & = \bX,
\end{align}
\end{subequations}
where $\bPcal_{\operatorname{\ell}} \in \C^{n \times n}$ and
$\bPcal_{\operatorname{r}} \in \C^{n \times n}$ are the left and right
projectors onto the subspace of finite eigenvalues of $\lambda \bE - \bA$.
In general, the matrices $\bPcal_{\operatorname{\ell}}$ and
$\bPcal_{\operatorname{r}}$ are given as spectral projectors via the
Weierstrass canonical form of $\lambda E - A$; see, for example,~\cite{MehS05}.
Typically, the computations necessary to obtain these projectors via subspace
decompositions are undesired in the case of large-scale sparse coefficient
matrices.
However, for several practically occurring matrix structures, the projectors
have been formulated explicitly in terms of parts of the coefficient matrices~\cite{BenS14, Sty08}.

In practice, the implicit application of $\bPcal_{\operatorname{\ell}}$ and
$\bPcal_{\operatorname{r}}$ using structural projectors is preferred over
the explicit formation and application of the projectors.
In this case, the solution of~\cref{eqn:projLyap} is directly computed on the
correct lower dimensional subspace.
Similar to the use of the spectral projectors, the implicit projection can only
be realized for certain matrix structures, for which the projectors onto the
correct subspaces and truncation of the coefficient matrices are known by
construction; see, for example,~\cite{BaeBSetal15, FreRM08, SaaV18}.
We refer the reader to the \matlab{} implementation of the Matrix Equations,
Sparse Solvers (\mbox{M-M.E.S.S.}) library~\cite{BenKS21, SaaKB23} for details
on the implementation of these structural projectors.

In all cases, we can note that none of the fundamental steps of the algorithms
presented in this paper change in the case of non-invertible $\bE$ matrices.
This case can be implemented by simply modifying the matrix-matrix and
matrix-vector operations needed in the algorithms to work directly on the
correct subspaces.
We also note that the invertibility of $\bE$ does not influence the role and
effect of the tangential directions in any way, which is why we continue with
the presentation of our results for the case of invertible $\bE$ matrices.


\section{Adaptive parameter selection}%
\label{sec:param}

As in the classical block ADI method (\Cref{sec:prelim}), the performance
of the tangential approach presented in \Cref{sec:tangADI} is governed by
the choice of iteration parameters, namely the shifts
$\alpha_{1}, \ldots, \alpha_{\ell} \in \C_{-}$ and the tangential directions
$\bt_{1}, \ldots, \bt_{\ell} \in \C^{m}$.
In the following, we present two adaptive frameworks that allow for the
selection of the shift parameters as well as the tangential directions during
the execution of the ADI method independent of each other.
In particular, the independence of the methods allows us to replace, for
example, the shift selection method by other approaches from the
literature~\cite{Wac13, Pen00a} as desired.


\subsection{Shift selection via projection}%
\label{sec:projshifts}

The selection of suitable shifts is a crucial point in any ADI-type method
and determines the practical convergence speed.
Poorly chosen shifts can result in many iteration steps, negating the
computational benefits of the low-rank formulation in the large-scale sparse
case. Classically, ADI theory relies on selecting shifts based on the
relative progression of the approximation error of the ADI iterates.
As detailed in~\cite{Kue16, Wac13}, minimizing the norm of this error
expression leads to the well-known ADI shift parameter problem, where optimal
shifts are given as the solution to the rational minimax problem
\begin{align} \label{eqn:minimax}
  \min\limits_{\{ \alpha_{1}, \dots, \alpha_{j} \}\subset \C_{-}}
    \left(\max\limits_{1 \leq i \leq n} \prod_{k = 1}^{j}
    \left\lvert \frac{\lambda_{i} - \overline{\alpha_{k}}}%
    {\lambda_{i} + \alpha_{k}} \right\rvert \right),
  \quad\text{where}\quad \lambda_{i} \in \Lambda(\bA, \bE).
\end{align}

Using~\cref{eqn:minimax} for the selection of shifts comes with two apparent
problems:
(i) the expression~\cref{eqn:minimax} is independent of the constant term in
the Lyapunov equation~\cref{eqn:lyap}, which is known to have a significant
influence on the numerical rank structure of the solution $\bX$, and
(ii) the solution to~\cref{eqn:minimax} can only be determined effectively in
very few cases of matrix spectra~\cite{Wac13}.
To mitigate these issues, the authors of~\cite{BenKS14} proposed an efficient
adaptive shift generation scheme for the classical block ADI approach
(\Cref{sec:prelim}) with $\bR = \eye{n}$ that has been observed to
provide suitable ADI shift parameters for general spectra of
$\lambda \bE - \bA$.
Thereby, the shifts are generated as Ritz values of the matrix pencil
$\lambda \bE - \bA$ projected onto the subspace spanned by the last ADI update
matrix~$\bV_{j - 1}$.
An extension to this was proposed in~\cite{BenSU16}, which suggests to further
enrich the projection space by including a set number $k$ of previous update
blocks $\bV_{j - k}, \bV_{j - k + 1}, \ldots, \bV_{j - 1}$ and then select a
subset of the computed Ritz values as shifts using a simplified variant
of~\cref{eqn:minimax}.
In contrast to the works~\cite{BenKS14, BenSU16}, we consider here the case
of Lyapunov equations with indefinite constant terms.
However, the projection shift method does not have to change for this case
since it is based on the range of the ADI updates, which is still determined
by the update matrices; cf. \Cref{sec:prelim}.

For the tangential ADI approach (\Cref{sec:tangADI}), the original
projection shifts generation from~\cite{BenKS14} would be insufficient.
Since the update terms are only single columns, a new shift would need to be
generated in every ADI step, which becomes computationally demanding, and
since only single shifts are generated, there would be no complex conjugate
shifts in the case of real coefficient matrices~\cref{eqn:lyap_real}.
The latter issue was the main motivation for the extension of the approach
in~\cite{BenSU16}, where real shifts alone did not lead to any meaningful
convergence of the ADI method due to the mechanical matrix structures.
Therefore, we also recommend in this work the use of the $k$ last update
columns $\bv_{j - k}, \bv_{j - k + 1}, \ldots, \bv_{j - 1}$ for the construction
of the projection space with some suitable $k \geq 1$.
Let $\bU \in \C^{n \times k}$ be a unitary basis of the subspace
$\mspan(\bv_{j - k}, \bv_{j - k + 1}, \ldots, \bv_{j - 1})$, then the projection
shifts are computed as the eigenvalues of the projected matrix pencil
$\lambda \bU^{\herm} \bE \bU - \bU^{\herm} \bA \bU$.
Note that depending on the coefficient matrices $\bA$ and $\bE$, some of the
projected eigenvalues may not be suited as ADI shifts.
Undesired shifts that lie in the right closed half-plane or are
infinite should be removed from the set.
Also, in principle, one may want to use some additional vectors to further
enrich the projection space beyond the number of actually desired shifts.
Then, similar to~\cite{BenSU16}, we propose to apply a modified version of
the minimax problem~\cref{eqn:minimax} to determine, which subset of shifts
should be used.
The resulting projection shift routine is outlined in \Cref{alg:projshifts}.
In the case that too few ADI update vectors are given or that a significant
amount of projected eigenvalues are invalid shifts, \Cref{alg:projshifts_else} in \Cref{alg:projshifts} is used
to return the remaining valid shifts.

\begin{algorithm}[t]
  \SetAlgoHangIndent{1pt}
  \DontPrintSemicolon

  \caption{Projection shift routine.}%
  \label{alg:projshifts}

  \KwIn{ADI update vectors $\bv_{j - k}, \bv_{j - k + 1}, \ldots, \bv_{j - 1}$,
    coefficient matrices $\bA, \bE$,
    number of desired shifts $\ell$.}
  \KwOut{Projection shifts $\alpha_{1}, \ldots, \alpha_{\ell}$.}

  Compute a unitary basis $\bU \in \C^{n \times k}$ of
    $\mspan(\bv_{j - k}, \bv_{j - k + 1}, \ldots, \bv_{j - 1})$.\;

  Compute the eigenvalues $\hat{\lambda}_{1}, \ldots, \hat{\lambda}_{k}$ of the
    matrix pencil $\lambda \bU^{\herm} \bE \bU - \bU^{\herm} \bA \bU$.

  Remove unwanted eigenvalues via
    $\{ \hat{\lambda}_{1}, \ldots, \hat{\lambda}_{p} \} =
    \{ \hat{\lambda}_{1}, \ldots, \hat{\lambda}_{k} \} \cap \C_{-}$.\;

  \eIf{$p > \ell$}{
    Select the $\ell$ most important shifts as the solution to
      \vspace{-.5\baselineskip}
      \begin{align*}
         \{ \alpha_{1}, \ldots, \alpha_{\ell} \} = 
           \argmin\limits_{\{ \alpha_{1}, \dots, \alpha_{\ell} \} \subset
             \{ \hat{\lambda}_{1}, \ldots, \hat{\lambda}_{p} \}}
             \left( \max\limits_{1 \leq i \leq p} \prod_{j = 1}^{\ell}
             \left\lvert \frac{\hat{\lambda}_{i} - \overline{\alpha_{j}}}%
             {\hat{\lambda}_{i} + \alpha_{j}} \right\rvert \right).
      \end{align*}\;
      \vspace{-1.5\baselineskip}
  }{
    Set $\ell \gets p$ and
      $\alpha_{i} \gets \hat{\lambda}_{i}$, for $i = 1, \ldots, \ell$.\;
    \label{alg:projshifts_else}
  }
\end{algorithm}

The remaining open problem for using the projection shifts as described in \Cref{alg:projshifts} is the initial phase of the ADI iteration when
no ADI update vectors $\bv_{j - k}, \bv_{j - k + 1}, \ldots, \bv_{j - 1}$ are
available yet.
One remedy to this problem is the use of alternative shift generation techniques
such as heuristic shifts~\cite{Pen00a} until enough ADI update vectors have
been generated to apply the projection shifts.
On the contrary, it has been suggested in~\cite{BenKS14} to use the
constant term of the Lyapunov equation~\cref{eqn:lyap} to generate the initial
projection space.
While this may work well as long as $m \leq \ell$, where $\ell$ is the number
of desired projection shifts, we consider cases in this work when $m > \ell$
and the rank of the constant term becomes too large to use its entirety for the
construction of the projection space.
In these situations, we propose to use some type of low-rank approximation
of the constant term in~\cref{eqn:lyap} to construct the initial projection
basis, for example, via the sparse singular value decomposition or randomized
sketching~\cite{HalMT11}.


\subsection{Tangential direction selection}%
\label{sec:tangdirect}

In addition to the ADI shifts, tangential directions have to be selected for
the presented tangential ADI method (\Cref{sec:tangADI}).
To this end, we first outline general restrictions on the choice of potential
tangential directions and derive specialized versions of \Cref{alg:LDLTADIcplx,alg:LDLTADIreal} based on these
results before we propose a new adaptive mechanism for efficiently selecting
suitable directions in the tangential ADI method.


\subsubsection{Challenges of static tangential directions}%
\label{sec:tangChange}

While both types of iteration parameters, ADI shifts and tangential directions,
appear in a similar way in \Cref{alg:LDLTADIcplx,alg:LDLTADIreal}, their influence on the convergence of the ADI method
is drastically different.
Most notably, the ADI method in any form is known to converge to the desired
solution even when only a single shift is used, i.e., when
$\alpha = \alpha_{1} = \ldots = \alpha_{\ell}$;
see, for example,~\cite{Pen00a, Wac13}.
The use of different shifts during the iteration only serves the purpose of
increasing the convergence speed, which is of particular importance when the
method is used to construct low-rank approximations.

On the contrary, using a variety of tangential directions is necessary in most
situations to allow convergence of the tangential ADI method towards the
correct solution.
Let us assume the extreme case of choosing all tangential directions to be
identical $\bt = \bt_{1} = \ldots = \bt_{\ell}$. The tangential ADI computes the solution approximation via the
summation
\begin{equation} \label{eqn:sol_tang}
  \bX_{j} = \sum\limits_{k = 1}^{j}
    \frac{-2 \real(\alpha_{k})}{\bt^{\herm} \bR^{-1} \bt}
    \bv_{k} \bv_{k}^{\herm},
\end{equation}
with the ADI update vectors $\bv_{k}$.
This formula is algebraically identical to the classical block ADI solution
update from \Cref{sec:prelim}, given by
\begin{equation} \label{eqn:sol_block}
  \bX_{j} = \sum\limits_{k = 1}^{j} -2\real(\alpha_{k}) \bV_{k}
    \tilde{\bR} \bV_{k}^{\herm},
\end{equation}
with the update vectors $\bV_{k} = \bv_{k}$ and the constant scalar center
term $\widetilde{\bR} = (\bt^{\herm} \bR^{-1} \bt)^{-1}$.
Consequently, applying the tangential ADI to~\cref{eqn:lyap} with a fixed
direction $\bt$ is mathematically equivalent to using the classical block
ADI method to solve instead
\begin{equation}\label{eqn:lyap_degen}
  \bA\bX\bE^{\herm} + \bE\bX\bA^\herm +
    \bB \bt ( \bt^{\herm} \bR^{-1} \bt )^{-1} \bt^{\herm} \bB^{\herm} = 0.
\end{equation}%
The solutions to~\cref{eqn:lyap} and~\cref{eqn:lyap_degen} are only identical,
if the constant terms in both equations are identical, which can only be the
case if the original constant term in~\cref{eqn:lyap} has rank $1$ and the tangential direction is chosen to be an eigenvector of~$\bR$.
Thus, in most scenarios, multiple tangential directions are needed for the
tangential ADI to converge to the desired solution.

Furthermore, the preemptive generation of suitable tangential directions similar
to the idea of heuristic shifts~\cite{Pen00a} is a highly difficult task
since it is not known in which way the tangential directions contribute
to the solution factors constructed by the ADI in terms of the approximation
error.
In the case that the constant term in~\cref{eqn:lyap} has a steep singular
value decay, it can be effectively approximated via an indefinite low-rank
factorization of the form
\begin{equation} \label{eqn:lowrankprod}
\widehat{\bB} \widehat{\bR} \widehat{\bB}^{\herm} = \sum\limits_{k = 1}^{r}
\bB \bt_{k} (\bt_{k}^{\herm} \bR^{-1} \bt_{k})^{-1} \bt_{k}^{\herm} \bB^{\herm} \approx
\bB \bR \bB^{\herm},
\end{equation}%
where $\widehat{\bB} \in \C^{n \times r}$ and
$\widehat{\bR} \in \C^{r \times r}$ with $r \ll m$.
Using the tangential directions in~\cref{eqn:lowrankprod} in a cyclic manner
is expected to allow the tangential ADI to converge as desired since the
principal information of the constant term is captured.
However, when $r$ is in a similar order of magnitude as $m$, this will lead
to similar computational costs and convergence behavior as the classical block
ADI approach.

In general, good choices of tangential directions that allow for a fast
convergence of the ADI method towards the correct solution are expected to
depend on the chosen ADI shifts as well as on the corresponding ADI iterates.
Therefore, an adaptive approach is needed to generate tangential directions
during the ADI iteration that take these factors into account.


\subsubsection{Eigenvector directions}%
\label{sec:tangEig}

Tangential directions as introduced in \Cref{sec:tangADI} compress the block
updates from the original ADI in \Cref{sec:prelim} down to single columns
and scalars.
Following the solution representations in~\cref{eqn:sol_tang}
and~\cref{eqn:sol_block}, the compression during the ADI iteration can be
simplified expressed in the form
\begin{equation} \label{eqn:approxfactor}
  \sum_{i = 1}^{j} \bV_{i} \bR \bV_{i}^{\herm} \approx
    \sum\limits_{k = 1}^{p} \bV_{k} \bt_{k} (\bt_{k}^{\herm} \bR^{-1} \bt_{k})^{-1}
    \bt_{k}^{\herm} \bV_{k}^{\herm}.
\end{equation}
At the center of this approximation stand the inverse
rank-$1$ projections of the indefinite center term~$\bR$.
If we assume for simplicity constant outer factors, the compression via
the tangential directions becomes essentially \begin{equation*}
  \bV_{i} \bR \bV_{i}^{\herm} \approx \sum\limits_{k = 1}^{p} \bV_{i} \bt_{k}
    (\bt_{k}^{\herm} \bR^{-1} \bt_{k})^{-1} \bt_{k}^{\herm} \bV_{i}^{\herm}
    = \bV_{i} \left( \sum\limits_{k = 1}^{p} \bt_{k}
    (\bt_{k}^{\herm} \bR^{-1} \bt_{k})^{-1} \bt_{k}^{\herm} \right)
    \bV_{i}^{\herm}.
\end{equation*}%
This suggests that tangential directions should be chosen in such a way that
their projections effectively decompose the $\bR$ matrix into rank-$1$ terms,
which in sum yield a good approximation.
In particular, we see that random directions are theoretically not expected to
yield any good results.

Since $\bR \in \C^{m \times m}$ is assumed to be Hermitian
(or symmetric in the case of real coefficients), the optimal rank-$1$
decomposition of the matrix with exactly $m$ terms is given by the
eigendecomposition of $\bR$.
Let $\bR = \bT \bS \bT^{\herm}$ be the eigendecomposition of $\bR$,
where $\bT = \begin{bmatrix} \bt_{1} & \ldots & \bt_{m} \end{bmatrix}$ is
the unitary matrix of eigenvectors and $\bS = \mdiag(s_{1}, \ldots, s_{m})$ is
the diagonal matrix of eigenvalues.
Then, we have for the eigenvectors and eigenvalues that
$\bt_{k}^{\herm} \bt_{k} = 1$ and $s_{k} = \bt_{k}^{\herm} \bR \bt_{k}$ so that
it holds
\begin{equation} \label{eqn:Rdecomp}
  \bR = \sum\limits_{k = 1}^{m} s_{k} \bt_{k} \bt_{k}^{\herm} =
    \sum\limits_{k = 1}^{m} \bt_{k} (\bt_{k}^{\herm} \bR \bt_{k}) \bt_{k}^{\herm}
  {}={} \sum\limits_{k = 1}^{m} \bt_{k} (\bt_{k}^{\herm}
    \bR^{-1} \bt_{k})^{-1} \bt_{k}^{\herm}.
\end{equation}
Inserting the eigenvectors into the general
approximation~\cref{eqn:approxfactor} then yields best low-rank approximations
of the terms occurring in the ADI iteration.

\begin{remark}[Reconstruction of block ADI updates]
  The original block ADI method (\Cref{sec:prelim}) can be recovered from
  the tangential approach (\Cref{sec:tangADI}) by making use of the
  eigendecomposition~\cref{eqn:Rdecomp}.
  To this end, for each ADI shift, we cycle through all $m$ eigenvectors
  of $\bR$ as tangential directions so that we recover update blocks of
  size~$m$ as in the classical block ADI.
  Note that the remaining difference between tangential and block ADI
  in the constructed matrices stems from the unitary transformation that
  diagonalizes the center matrix $\bD_{j}$ in the tangential approach,
  while in the block method that center matrix is block diagonal.
\end{remark}

The algorithms presented in \Cref{sec:tangADI} were formulated for the
general case of tangential directions.
With the restriction of these directions to the eigenvectors of $\bR$, these
algorithms can be further adapted.
Due to the practical importance of the choice of these specific tangential
directions and the corresponding algorithmic changes, we outline the tangential
ADI method with eigenvector directions for the cases of complex and real
coefficient matrices in \Cref{alg:LDLEigTADIcplx,alg:LDLEigTADIreal}.
In comparison to the general \Cref{alg:LDLTADIcplx,alg:LDLTADIreal},
the eigenvector-based variants in \Cref{alg:LDLEigTADIcplx,alg:LDLEigTADIreal}
fully avoid the rescaling with the inverse center matrix $\bR^{-1}$ and
instead involve the eigenvalues of $\bR$.
Also, the selection of tangential directions in
\Cref{alg:LDLEigTADIcplx_tangselect} of \Cref{alg:LDLEigTADIcplx}
and \Cref{alg:LDLEigTADIreal_tangselect} of \Cref{alg:LDLEigTADIreal}
simplifies to the discrete choice between $m$ different fixed vectors.
We will make use of this in the next section to devise multiple schemes that
allow the adaptive selection of tangential directions in every ADI step.

\begin{algorithm}[t]
  \SetAlgoHangIndent{1pt}
  \DontPrintSemicolon

  \caption{Indefinite tangential ADI with eigenvectors: complex case.}
  \label{alg:LDLEigTADIcplx}

  \KwIn{Coefficient matrices $\bA, \bE, \bB, \bR$ from~\cref{eqn:lyap},
    shift parameters $\alpha_{1}, \ldots, \alpha_{\ell} \in \C_{-}$,
    residual tolerance $\epsilon > 0$.}
  \KwOut{Approximate solution factors $\bL_{j-1}, \bD_{j-1}$ so that
    $\bL_{j-1} \bD_{j-1} \bL_{j-1}^{\herm} \approx \bX$.}

  Initialize $\bW_{0} = \bB$, $\bL_{0} = [~]$, and $j = 1$.\;
  
  Compute the eigendecomposition $\bR = \bT \bS \bT^{\herm}$, with
    $\bT = \begin{bmatrix} \bt_{1} & \ldots & \bt_{m} \end{bmatrix}$ and
    $\bS = \mdiag(s_{1}, \ldots, s_{m})$.\;
  
  \While{$\lVert \bW_{j - 1} \bR \bW_{j - 1}^{\herm} \rVert \geq
    \epsilon \lVert \bB \bR \bB^{\herm} \rVert$}{
    Select eigenvector $\bt_{p_{j}}$ with corresponding eigenvalue
      $s_{p_{j}}$ for $p_{j} \in \{1, \ldots, m\}$.\;
    \label{alg:LDLEigTADIcplx_tangselect}
      
    Solve $(\bA + \alpha_{j} \bE) \bv_{j} = \bW_{j - 1} \bt_{p_{j}}$
      for $\bv_{j}$.\;

    Update the residual factor
      $\bW_{j} = \bW_{j - 1} - 2 \real(\alpha_{j})\bE \bv_{j}
      \bt_{p_{j}}^{\herm}$.\;

    Update the solution factor
      $\bL_{j} = \begin{bmatrix} \bL_{j - 1} & \bv_{j} \end{bmatrix}$.
    
    Increment $j \gets j + 1$.\;
  }

  Compute the solution factor
    $\bD_{j - 1} = -2 \mdiag(\real(\alpha_{1}) s_{p_{1}}, \ldots,
      \real(\alpha_{j - 1}) s_{p_{j - 1}})$.\;
\end{algorithm}

\begin{algorithm}[t]
  \SetAlgoHangIndent{1pt}
  \DontPrintSemicolon
  
  \caption{Indefinite tangential ADI with eigenvectors: real case.}
  \label{alg:LDLEigTADIreal}

  \KwIn{Real coefficient matrices $\bA, \bE, \bB, \bR$
    from~\cref{eqn:lyap_real},
    shift parameters $\alpha_{1}, \ldots, \alpha_{\ell} \in \C_{-}$ closed
    under conjugation and ordered as pairs,
    residual tolerance $\epsilon > 0$.}
  \KwOut{Approximate solution factors $\bL_{j - 1}, \bD_{j - 1}$ so that
    $\bL_{j - 1} \bD_{j - 1} \bL_{j - 1}^{\trans} \approx \bX$.}

  Initialize $\bW_{0} = \bB$, $\bL_{0} = [~]$, $\eta = \sqrt{2}$, and $j = 1$.\;
  
  Compute the eigendecomposition $\bR = \bT \bS \bT^{\trans}$, with
    $\bT = \begin{bmatrix} \bt_{1} & \ldots & \bt_{m} \end{bmatrix}$ and
    $\bS = \mdiag(s_{1}, \ldots, s_{m})$.\;

  \While{$\lVert \bW_{j - 1} \bR \bW_{j - 1}^{\herm} \rVert \geq
    \epsilon \lVert \bB \bR \bB^{\herm} \rVert$}{

    Select eigenvector $\bt_{p_{j}}$ with corresponding eigenvalue
      $s_{p_{j}}$ for $p_{j} \in \{1, \ldots, m\}$.\;
    \label{alg:LDLEigTADIreal_tangselect}
    
    Solve $(\bA + \alpha_{j} \bE) \bv_{j} = \bW_{j - 1} \bt_{p_{j}}$
      for $\bv_{j}$.\;

    \eIf{$\alpha_{j}$ is real}{
      Update the residual and the solution factors via
        \vspace{-.5\baselineskip}
        \begin{align*}
          \bW_{j} & = \bW_{j - 1} - 2 \alpha_{j} \bE \bv_{j}
            \bt_{p_{j}}^{\trans}, \\
          \bL_{j} & = \begin{bmatrix} \bL_{j - 1} & \bv_{j} \end{bmatrix}.
        \end{align*}\;
        \vspace{-1.5\baselineskip}
    }{
      Compute the scaling factor
        $\delta_{j} = \real(\alpha_{j}) / \imag(\alpha_{j})$.\;

      Update the residual and the solution factors via
        \vspace{-.5\baselineskip}
        \begin{align*}
          \bW_{j + 1} & = \bW_{j - 1} - 4 \real(\alpha_{j}) \bE \big(
            \real(\bv_{j}) + \delta_{j} \imag(\bv_{j}) \big)
            \bt_{p_{j}}^{\trans}, \\
          \bL_{j + 1} & = \begin{bmatrix} \bL_{j - 1} &
            \eta \big( \real(\bv_{j}) + \delta_{j} \imag(\bv_{j}) \big) &
            \eta \sqrt{\delta_{j}^{2} + 1} \imag(\bv_{j}) \end{bmatrix}.
        \end{align*}\;
        \vspace{-1.5\baselineskip}

      Increment $j \gets j + 1$.\;
    }

    Increment $j \gets j + 1$.\;
  }

  Compute the solution factor
    $\bD_{j - 1} = -2 \mdiag(\real(\alpha_{1}) s_{p_{1}}, \ldots,
    \real(\alpha_{j - 1}) s_{p_{j - 1}})$.\;
\end{algorithm}

Before moving on, we would like to outline the practical difference
between the standard literature approaches of selecting potentially arbitrary
tangential directions and our proposed strategy of restricting these directions to
the eigenvectors of the center matrix $\bR$.
For this purpose, we are using a randomly generated stable test example with
outer dimension $n = 1\,000$ and with a constant term of rank $m = 10$
that has an indefinite $\bR$ matrix.
We compare the classical block ADI from \Cref{sec:prelim} with the
tangential approaches using random sets of orthogonalized directions and the
eigenvectors of $\bR$.
The convergence behavior of the different methods over the size of the solution
factors is shown in \Cref{fig:tangchoice}.
The classical block ADI (\adi{}) converges to a reasonably small normalized
Lyapunov residual computed as
$\lVert \bLcal_{j} \rVert_{2} / \lVert \bB \bR \bB^{\herm} \rVert_{2}$
showcasing that the ADI can be successfully applied here.
The tangential ADI with eigenvectors (\tadieig{}) converges similarly fast but
the tangential ADI with random directions (\tadirand{}) diverges instead.
This fits our previous discussion that random directions are not suited
to solve the approximation problem~\cref{eqn:approxfactor}.

\begin{figure}[t]
  \centering
  \tikzexternalenable%
  \tikzsetnextfilename{no_separation}%
  \begin{tikzpicture}[font = \plotfontsize\normalfont]

  \pgfplotstableread{graphics/data/exp1_efficiency_adi.dat}\tableADI
  \pgfplotstableread{graphics/data/exp1_efficiency_eigtadi.dat}\tableEigTADI
  \pgfplotstableread{graphics/data/exp1_efficiency_rand.dat}\tableRandTADI

  \begin{semilogyaxis}[%
  scale only axis,
  width               = .8\linewidth,
  height              = .2\linewidth,
  xmin                = 0,
  xmax                = 300,
  ymin                = 1e-12,
  ymax                = 1e+4,
  xtick               = {0, 50, 100, 150, 200, 250, 300},
  ytick               = {1e+4, 1e+0, 1e-4, 1e-8, 1e-12},
  xlabel              = {dimension of solution factor $\bL_{j}$},
  ylabel              = {normalized residual},
  ylabel style        = {yshift = -.3em},
  scaled x ticks      = false,
  x tick label style  = {/pgf/number format/1000 sep={\,}},
  y tick label style  = {/pgf/number format/1000 sep={\,}},
  grid                = major,
  legend cell align   = left,
  legend style        = {
    at        = {(1, 1)},
    anchor    = north east,
    align     = left}
  ]

    \addplot[adi_imp_style] 
      table[x = x_iter_rank, y = y_implicit] {\tableADI};
    \addlegendentry{\adi{}}

    \addplot[tadi_rand_style, mark repeat = 3, mark phase = 0] 
      table[x = x_iter_rank, y = y_implicit] {\tableRandTADI};
    \addlegendentry{\tadirand{}}
  
    \addplot[eigtadi_imp_style, mark repeat = 25, mark phase = 5] 
      table[x = x_iter_rank, y = y_implicit] {\tableEigTADI};
    \addlegendentry{\tadieig{}}  
  \end{semilogyaxis}
\end{tikzpicture}%
  \tikzexternaldisable%

  \caption{Convergence behavior of the classical block ADI method (\adi{}), the
    tangential ADI with random directions (\tadirand{}) and
    the tangential ADI with eigenvector directions (\tadieig{}):
    The eigenvector-based method behaves similarly to the classical approach
    and converges to a small normalized residual, while the choice of random
    directions leads to divergence of the tangential ADI.}
  \label{fig:tangchoice}
\end{figure}

We note that throughout various numerical experiments, we have observed that
in the case of indefinite $\bR$ matrices, the tangential ADI always diverges as
shown in \Cref{fig:tangchoice} when the tangential directions are chosen as
linear combinations of eigenvectors corresponding to at least one positive and
one negative eigenvalue.
Random tangential directions satisfy this property with high probability when
$\bR$ is indefinite.
On the other side, if the tangential directions have been chosen as arbitrary
vectors from either the eigenspace corresponding to the positive eigenvalues or
the eigenspace corresponding to the negative eigenvalues, the convergence
of the tangential ADI essentially stagnates about the initial residual.
Overall, we could not observe any choice of tangential directions other than
the eigenvectors of $\bR$ that resulted in visible convergence of the
tangential ADI method.
Finally, the special case of trivial center matrices is discussed in the
following remark.

\begin{remark}[Relation to the previous tangential ADI]
  The previously proposed tangential ADI approach from~\cite{WolPL16} can be
  recovered from this work by setting $\bR = \eye{m}$ in the Lyapunov
  equation~\cref{eqn:lyap_real} and by normalizing all tangential directions
  $\bt_{j} \in \R^{m}$ so that $\bt_{j}^{\trans} \bt_{j} = 1$.
  Because any non-zero vector is an eigenvector of the identity matrix,
  restricting the directions to standard canonical unit vectors is unnecessary and
  would overly restrict the tangential updates to picking single columns of the
  ADI block factor.
  Instead, the approximation problem~\cref{eqn:approxfactor} simplifies to
  \begin{equation*}
    \sum\limits_{i = 1}^{j} \bV_{i} \bV_{i}^{\herm} \approx
    \sum\limits_{k = 1}^{p} \bV_{k} \bt_{k} \bt_{k}^{\herm} \bV_{k}^{\herm},
  \end{equation*}
  for which any choice of normalized tangential directions will yield reasonable
  results, as illustrated in~\cite{WolPL16}.
  Furthermore, if the original center matrix $\bR$ is strictly positive definite
  rather than indefinite, one can simply absorb its Cholesky factor into the
  outer matrix $\bB$ to recover the $\bR = \eye{m}$ case, allowing the use of
  unrestricted directions without potential convergence issues mentioned before.
  We note that the normalization used in~\cite{WolPL16} is not necessary in our
  work since the corresponding rescaling is integrated into
  \Cref{alg:LDLTADIcplx,alg:LDLTADIreal}.
  Additionally, this discussion holds when using any trivially scaled variant
  of the identity for the center term, i.e., when
  $\bR = \beta \eye{m}$ for some~$\beta \in \C$.
\end{remark}


\subsubsection{Adaptive selections via heuristic scoring}%
\label{sec:tangScores}

With the tangential directions restricted to the eigenvectors of $\bR$,
the problem of adaptively constructing directions during the ADI iteration
simplifies to the selection of $1$ out of $m$ pre-computed vectors in every
tangential ADI step.
In general, the goal is to select the tangential direction $\bt_{j}$ from the
eigenvectors of $\bR$, which for a given shift $\alpha_{j}$ and a given
previous residual factor $\bW_{j - 1}$ minimizes the approximation error
of the constructed solution vectors via the resulting update vector $\bv_{j}$,
which is given by $\bv_{j} = (\bA + \alpha_{j} \bE)^{-1} \bW_{j - 1} \bt_{j}$.
The exact solution to this problem is practically infeasible, since it requires
the construction of all $m$ possible update vectors and the computation
of the corresponding residuals, which is computationally much more costly than
performing classical block ADI steps.
Instead, we propose the use of heuristic values that indicate the importance
of potential updates so that we may base the choice of tangential directions
on this heuristic.

In general, we expect all update vectors $\bv_{j}$ computed via the
eigenvectors of $\bR$ to positively contribute towards the solution of the
Lyapunov equation in some way since all those update vectors are constructed
following the classical ADI theory.
To quantify the contribution, we may simply use the Euclidean norm of the
update vectors since a large vector norm indicates large changes in the
solution approximation while a small vector norm indicates only minor changes.
Thus, a heuristic approach for choosing the tangential directions is given
by the solution to the problem
\begin{equation} \label{eqn:fullheur}
  \bt_{j} = \argmax\limits_{\bt \in \{\bt_{1}, \ldots, \bt_{m}\}}
    \lVert \bv_{j} \rVert_{2} = 
    \argmax\limits_{\bt \in \{\bt_{1}, \ldots, \bt_{m}\}}
    \lVert (\bA + \alpha_{j} \bE)^{-1} \bW_{j - 1} \bt \rVert_{2},
\end{equation}
where $\bt_{1}, \ldots, \bt_{m}$ are the unitary eigenvectors of $\bR$.

Similar to the original problem of selecting the optimal tangential direction,
this heuristic~\cref{eqn:fullheur} is practically infeasible as it requires the
solution of $m$ linear systems of dimension $n$.
To lower the computational costs of the heuristic, we propose the use of the
same Galerkin projection framework that is used in the generation of
adaptive shifts in \Cref{sec:projshifts}.
To this end, let $\bU \in \C^{n \times k}$ denote the unitary basis
matrix for the $k$-dimensional subspace spanned by the previous update
vectors $\bv_{j - k}, \bv_{j - k + 1}, \ldots, \bv_{j - 1}$.
Then, compute the next tangential direction as the solution to the problem
\begin{equation} \label{eqn:romheur}
  \bt_{j} = \argmax\limits_{\bt \in \{\bt_{1}, \ldots, \bt_{m}\}}
    \lVert (\widehat{\bA} + \alpha_{j} \widehat{\bE})^{-1}
    \widehat{\bW}_{j - 1} \bt \rVert_{2},
\end{equation}
where $\widehat{\bA} = \bU^{\herm} \bA \bU$,
$\widehat{\bE} = \bU^{\herm} \bE \bU$,
$\widehat{\bW}_{j - 1} = \bU^{\herm} \bW_{j - 1}$, and
where $\bt_{1}, \ldots, \bt_{m}$ are the unitary eigenvectors of $\bR$.
This projective heuristic~\cref{eqn:romheur} requires only the solution to
$m$ linear systems of dimension $k$, where $k$ is typically chosen to be
very small in comparison to the full dimension $n$.
Note that these $m$ linear systems only differ in their right-hand side so that
all $m$ solutions can be computed simultaneously using a single LU
decomposition.
As in the case of the projection shifts (\Cref{sec:projshifts}), we expect the
last ADI update columns to contain enough information about the iteration
behavior to allow for further improvements.
Due to its importance for practical implementations, we have summarized the
algorithmic description of the projection-based heuristic selection of the
tangential directions in \Cref{alg:projtang}.

\begin{algorithm}[t]
  \SetAlgoHangIndent{1pt}
  \DontPrintSemicolon

  \caption{Projection-based heuristic tangential direction selection.}
  \label{alg:projtang}

  \KwIn{Unitary eigenvectors $\bt_{1}, \ldots, \bt_{m}$ of $\bR$,
    previous ADI update vectors
    $\bv_{j - k}, \bv_{j - k + 1}, \ldots, \bv_{j - 1}$,
    previous residual vector $\bW_{j - 1}$,
    next ADI shift $\alpha_{j}$,
    coefficient matrices $\bA, \bE$.}
  \KwOut{Next tangential direction $\bt$.}

  Compute a unitary basis matrix $\bU \in \C^{n \times k}$ of
    $\mspan(\bv_{j - k}, \bv_{j - k + 1}, \ldots, \bv_{j - 1})$.\;

  Compute the projected terms $\widehat{\bA} = \bU^{\herm} \bA \bU$,
    $\widehat{\bE} = \bU^{\herm} \bE \bU$, and
    $\widehat{\bW}_{j - 1} = \bU^{\herm} \bW_{j - 1}$.

  Compute the projected update vectors
    $\widehat{\bV} = \begin{bmatrix} \hat{\bv}_{1} & \ldots & \hat{\bv}_{m}
    \end{bmatrix}$ via the solution of
    \vspace{-.5\baselineskip}
    \begin{equation*}
      (\widehat{\bA} + \alpha_{j} \widehat{\bE}) \widehat{\bV} = 
        \widehat{\bW}_{j - 1} \begin{bmatrix} \bt_{1} & \ldots & \bt_{m}
        \end{bmatrix}.
    \end{equation*}\;
    \vspace{-1.5\baselineskip}

  Select the index with the largest vector norm
    $p = \argmax_{\ell = 1, \ldots, m} \lVert \hat{\bv}_{\ell} \rVert_{2}$.\;

  Return the corresponding tangential direction $\bt \gets \bt_{p}$.\;
\end{algorithm}

\begin{remark}[Simplified heuristics]
  The motivation for the heuristics in~\cref{eqn:fullheur} and~\cref{eqn:romheur} was to
  determine the tangential directions that contribute the most towards the
  solution approximation.
  In contrast to that, one may consider to use instead the tangential directions
  that capture most of the current residual factor to drive the residual
  towards zero.
  A corresponding heuristic for selecting tangential directions is given by
  \begin{equation} \label{eqn:resheur}
    \bt_{j} = \argmax\limits_{\bt \in \{ \bt_{1}, \ldots, \bt_{m} \}}
      \lVert \bW_{j - 1} \bt \rVert_{2},
  \end{equation}
  where $\bW_{j - 1}$ is the previous residual factor and
  $\bt_{1}, \ldots, \bt_{m}$ are the unitary eigenvectors of $\bR$.
  This heuristic is computationally inexpensive even on the high-dimensional
  residual factor since it fully avoids the solution of linear systems.
  However, Equation~\cref{eqn:resheur} does not account for the influence of the next
  ADI shift so that we expect a tangential ADI algorithm equipped
  with~\cref{eqn:resheur} to converge slower in terms of iteration steps than a
  tangential ADI that uses~\cref{eqn:romheur}.
\end{remark}

To illustrate the difference between the proposed selection heuristics,
we employ the same test example that has been used in \Cref{sec:tangEig}.
The results shown in \Cref{fig:heuristics} verify that the full heuristic
from~\cref{eqn:fullheur} performs consistently best with the smallest number
of iterations.
However, the projected variant~\cref{eqn:romheur}, which is computationally
significantly cheaper, performs very similarly to the full approach.
Finally, the residual-based variant from~\cref{eqn:resheur} performs the worst of all the considered strategies but still yields reasonably fast convergence
so that it could be employed in practice.
We also note that the residual-based approach shows several phases of
stagnation throughout the iteration.
This implies that selecting directions with respect to the largest residual
contribution does not necessarily lead to a reduction of the residual norm.

\begin{figure}[t]
  \centering
  \tikzexternalenable%
  \tikzsetnextfilename{scoring_methods}%
  \begin{tikzpicture}[font = \plotfontsize\normalfont]

\pgfplotstableread{graphics/data/exp2_tangent_selection_full.dat}\tableFULL
\pgfplotstableread{graphics/data/exp2_tangent_selection_proj.dat}\tablePROJ
\pgfplotstableread{graphics/data/exp2_tangent_selection_Wnorm.dat}\tableWNORM

\begin{semilogyaxis}[%
  scale only axis,
  width               = .8\linewidth,
  height              = .2\linewidth,
  xmin                = 0,
  xmax                = 350,
  ymin                = 1e-12,
  ymax                = 1e+1,
  xtick               = {0, 50, 100, 150, 200, 250, 300, 350},
  ytick               = {1e+0, 1e-4, 1e-8, 1e-12},
  xlabel              = {tangential ADI iteration step $j$},
  ylabel              = {normalized residual},
  ylabel style        = {yshift = -.3em},
  scaled x ticks      = false,
  x tick label style  = {/pgf/number format/1000 sep={\,}},
  y tick label style  = {/pgf/number format/1000 sep={\,}},
  grid                = major,
  legend cell align   = left,
  legend style        = {
    at        = {(1, 1)},
    anchor    = north east,
    align     = left}
  ]

    \addplot[eigtadi_fullheur, mark repeat = 25, mark phase = 0] 
      table[x = x_iter_rank, y = y_implicit] {\tableFULL};
    \addlegendentry{Full heuristic}
  
    \addplot[eigtadi_imp_style, mark repeat = 25, mark phase = 10] 
      table[x=x_iter_rank, y=y_implicit] {\tablePROJ};
    \addlegendentry{Projected heuristic}

    \addplot[eigtadi_resheur, mark repeat = 25, mark phase = 20] 
      table[x=x_iter_rank, y=y_implicit] {\tableWNORM};
    \addlegendentry{Residual heuristic}
  
  \end{semilogyaxis}
\end{tikzpicture}%
  \tikzexternaldisable%

  \caption{Comparison of the different heuristic selection strategies for
    tangential direction selection:
    The full heuristic selection based on~\cref{eqn:fullheur} performs best
    closely followed by the projected variant~\cref{eqn:romheur}.
    The residual-based heuristic selection~\cref{eqn:resheur} performs slightly
    worse than the other heuristics but still provides reasonably fast
    convergence.}
  \label{fig:heuristics}
\end{figure}


\section{Numerical experiments}%
\label{sec:numerics}

In this section, we demonstrate the proposed tangential ADI method (further on
denoted by \tadi{}) in comparison to the classical block approach
(from here on as \adi{}) on several synthetic and real-world benchmark examples.
Following the discussions in \Cref{sec:param}, we restrict ourselves
for \tadi{} to the tangential ADI method with eigenvector directions as outlined in \Cref{alg:LDLEigTADIcplx,alg:LDLEigTADIreal} and using
the projected heuristic selection of the tangential directions from \Cref{alg:projtang}.
Also, for \adi{} and \tadi{}, shifts are adaptively computed using the
projection shift routine from \Cref{alg:projshifts}.
For the comparisons, we are computing the normalized residual norms as
\begin{equation*}
  \frac{\lVert \bLcal_{j} \rVert_{\fro}}{\lVert \bB \bR \bB^{\herm} \rVert_{\fro}} =
    \frac{\lVert \bA \bX_{j} \bE^{\herm} +  \bE \bX_{j} \bA^{\herm} +
    \bB \bR \bB^{\herm} \rVert_{\fro}}{\lVert \bB \bR \bB^{\herm} \rVert_{\fro}},
\end{equation*}
where $\bX_{j} = \bL_{j} \bD_{j} \bL_{j}^{\herm}$ is the current ADI iterate.
In all experiments, all shown methods are run to compute approximations
with a normalized residual norm of less than $10^{-12}$.
In addition, in the large-scale experiments, we also compare \tadi{} and
\adi{} against the ADI and EigTADI approaches that employ intermediate
truncations of the solution factors.
We denote by \aditrunc{} the ADI method, which uses the truncation in every
step of the algorithm and by \taditrunc{} an EigTADI variant, which uses the
truncation of the solution in every $k = 6$ or $k = 12$ steps for the
(bi)linear rail and butterfly examples, respectively.
With this, we aim to demonstrate an alternate approach for saving memory during
the computations in comparison to \tadi{}.
However, we note that the use of intermediate truncation-based compressions can
lead to loss of accuracy for the true ADI residual following the theory
about the inexact ADI method developed in~\cite{KueF20}.

The reported numerical experiments have been performed on a MacBook Air with
16\,GB of RAM and an Apple M2 processor running macOS Ventura version 13.4 with
MATLAB 25.1.0.2973910 (R2025a) Update 1.
The source codes, data, and results of the numerical experiments reported in
this section and in \Cref{sec:param} are available at~\cite{supSmiW26}.


\subsection{Verifying convergence and implicit residuals}%
\label{sec:resconv}

Before tackling high-dimensional real-world examples from the literature,
we will verify the correctness of our theoretical results with two medium-scale
synthetic numerical examples.
Both examples are randomly generated with the dimensions $n = 1\,000$ and
$m = 20$ and modified so that the eigenvalues of the matrix pencils
$\lambda \bE - \bA$ lie in the open left half-plane.
The first example has complex coefficient matrices, while the second one
has real coefficients.
For the real example, we note that the matrices $\bA$ and $\bE$ are
non-symmetric so that the shifts generated via \Cref{alg:projshifts} can
be both real or occur in complex conjugate pairs.
For both examples, the center matrix $\bR$ is constructed to be
indefinite.
The final coefficient matrices for both examples have been exported and saved
in the accompanying code package~\cite{supSmiW26} for potential reproduction
of the presented experiments.

The two goals of the numerical experiments with these examples are to
verify for the complex as well as the real versions of the indefinite factorized
\adi{} and \tadi{} that (i) the methods do converge to a reasonable
approximation of the solution, and (ii) that the formulas from \Cref{sec:prelim,sec:tangADI} for the implicitly updated
residual factors do match with the true Lyapunov residual.
The results of these experiments are shown in \Cref{fig:synthetic}.
We see that for both examples, the classical block ADI and the novel
tangential approach converge to a normalized residual of about $10^{-12}$.
Also, for both types of methods, the implicitly updated residual factors align
with the true Lyapunov residuals, which verifies the computational correctness
of our theoretical results from \Cref{sec:prelim,sec:tangADI}.

\begin{figure}[t]
  \centering
  \begin{subfigure}[b]{.49\linewidth}
    \centering
  \tikzexternalenable%
  \tikzsetnextfilename{synthetic_complex}%
  \begin{tikzpicture}[font = \plotfontsize\normalfont]
  \pgfplotstableread{graphics/data/complex_residual_adi.dat}\tableADI
  \pgfplotstableread{graphics/data/complex_residual_eigtadi.dat}\tableEigTADI
  
  \begin{semilogyaxis}[%
    scale only axis,
    width              = .7\linewidth,
    height             = .5\linewidth,
    xmin               = 0,
    xmax               = 210,
    ymin               = 1e-12,
    ymax               = 1e+1,
    xtick              = {0, 50, 100, 150, 200},
    ytick              = {1e+0, 1e-4, 1e-8, 1e-12},
    xlabel             = {dimension of solution factor $\bL_{j}$},
    ylabel             = {normalized residual},
    ylabel style       = {yshift = -.3em},
    scaled x ticks     = false,
    x tick label style = {/pgf/number format/1000 sep={\,}},
    y tick label style = {/pgf/number format/1000 sep={\,}},
    grid               = major
  ]

    \addplot[adi_imp_style] 
      table[x=x_iter_rank, y=y_implicit] {\tableADI};

    \addplot[adi_exp_style] 
      table[x=x_iter_rank, y=y_explicit] {\tableADI};

    \addplot[eigtadi_imp_style, mark repeat = 25, mark phase = 0] 
      table[x=x_iter_rank, y=y_implicit] {\tableEigTADI};

    \addplot[eigtadi_exp_style, mark repeat = 25, mark phase = 15] 
      table[x=x_iter_rank, y=y_explicit] {\tableEigTADI};
    
  \end{semilogyaxis}
\end{tikzpicture}%
  \tikzexternaldisable%

    \caption{complex example}
  \end{subfigure}%
  \hfill%
  \begin{subfigure}[b]{.49\linewidth}
    \centering
  \tikzexternalenable%
  \tikzsetnextfilename{synthetic_real}%
  \begin{tikzpicture}[font = \plotfontsize\normalfont]
  \pgfplotstableread{graphics/data/real_nonsymmetric_residual_adi.dat}\tableADI
  \pgfplotstableread{graphics/data/real_nonsymmetric_residual_eigtadi.dat}\tableEigTADI

  \begin{semilogyaxis}[%
    scale only axis,
    width              = .7\linewidth,
    height             = .5\linewidth,
    xmin               = 0,
    xmax               = 170,
    ymin               = 1e-12,
    ymax               = 1e+1,
    xtick              = {0, 40, 80, 120, 160},
    ytick              = {1e+0, 1e-4, 1e-8, 1e-12},
    xlabel             = {dimension of solution factor $\bL_{j}$},
    ylabel             = {normalized residual},
    ylabel style       = {yshift = -.3em},
    scaled x ticks     = false,
    x tick label style = {/pgf/number format/1000 sep={\,}},
    y tick label style = {/pgf/number format/1000 sep={\,}},
    grid               = major
  ]

    \addplot[adi_imp_style] 
      table[x=x_iter_rank, y=y_implicit] {\tableADI};

    \addplot[adi_exp_style] 
      table[x=x_iter_rank, y=y_explicit] {\tableADI};

    \addplot[eigtadi_imp_style, mark repeat = 25, mark phase = 0] 
      table[x=x_iter_rank, y=y_implicit] {\tableEigTADI};

    \addplot[eigtadi_exp_style, mark repeat = 25, mark phase = 15] 
      table[x=x_iter_rank, y=y_explicit] {\tableEigTADI};
  \end{semilogyaxis}
\end{tikzpicture}%
  \tikzexternaldisable%

    \caption{real example}
  \end{subfigure}

  \vspace{.5\baselineskip}
  \tikzexternalenable%
  \tikzsetnextfilename{synthetic_legend}%
  \begin{tikzpicture}[font = \plotfontsize\normalfont]
  \begin{axis}[%
    hide axis,
    scale only axis,
    width          = .8\linewidth,
    height         = .2\linewidth,
    xmin           = 0,
    xmax           = 1,
    ymin           = 0,
    ymax           = 1,
    legend columns = 2, 
    legend style   = {
      at     = {(0,0)},
      anchor = center,
      /tikz/every even column/.append style = {column sep = 0.5cm}}
  ]
    
    \addplot[adi_imp_style] coordinates {(0,0)};
    \addlegendentry{\adi{}, implicit residual}

    \addplot[adi_exp_style] coordinates {(0,0)};
    \addlegendentry{\adi{}, explicit residual}

    \addplot[eigtadi_imp_style] coordinates {(0,0)};
    \addlegendentry{\tadi{}, implicit residual}

    \addplot[eigtadi_exp_style] coordinates {(0,0)};
    \addlegendentry{\tadi{}, explicit residual}
  \end{axis}
\end{tikzpicture}%
  \tikzexternaldisable%

  \caption{Convergence plots for ADI methods on synthetic examples:
    The tangential ADI methods converge similarly to the classical block ADI
    methods in these examples to a normalized residual of about $10^{-12}$.
    For all shown methods, the implicit residual update formulas align with the
    true explicit Lyapunov residuals.}
  \label{fig:synthetic}
\end{figure}


\subsection{Examples with large-scale coefficient matrices}%
\label{sec:largescale}

In this section, we consider three known benchmark examples from the literature
and compare the proposed tangential ADI method to the classical block
approach.


\subsubsection{Simplified heat transfer in a steel bar}%
\label{sec:rail}

The first example is the classical steel rail data set
from~\cite{morwiki_steel, Saa03}.
This data set describes the cooling process of a steel bar via a linear system
of ordinary differential equations of the form
\begin{equation*}
  \widetilde{\bE} \dot{\bx}(t) = \widetilde{\bA} \bx(t) +
    \widetilde{\bB} \bu(t), \quad \by(t) = \widetilde{\bC} \bx(t), 
\end{equation*}
where $\widetilde{\bE}, \widetilde{\bA} \in \R^{n \times n}$,
$\widetilde{\bB} \in \R^{n \times m}$, and
$\widetilde{\bC} \in \R^{p \times n}$, with $n = 20\,209$, $m = 7$, and
$p = 6$.
Therein, the matrix $\widetilde{\bE}$ is symmetric positive definite while
$\widetilde{\bA}$ is symmetric negative definite.
Lyapunov equations are used in the construction of reduced-order models and
controllers for this dynamical system.
Here, we aim for the solution of the observability Lyapunov equations, which
is given as the solution of~\cref{eqn:lyap_real} by setting
$\bA = \widetilde{\bA}$, $\bE = \widetilde{\bE}$, and
$\bB = \widetilde{\bC}^{\trans}$.
For the center matrix $\bR \in \R^{p \times p}$, we consider an indefinite
symmetric matrix inspired by the setting of
$\Hinf$-control~\cite{BenHW22, MusG91}, where the
positive eigenvalues correspond to the performance of the system while the
negative eigenvalues describe how perturbations are observed.

The performance of \tadi{}, \adi{}, and their dynamically truncated
variants can be seen in \Cref{fig:rail}.
Comparing the different convergence plots, we see that the tangential
approaches \tadi{} and \taditrunc{} provide consistently smaller solution factors
throughout the iteration compared to their block counterparts \adi{}
and \aditrunc{} for the same level of accuracy.
The corresponding computation times are shown in \Cref{tab:rail_timing}.
In addition to the basic solvers, a post truncation is performed
to reduce the final solution factors to a similar numerical rank for all
methods, with the corresponding computation times shown in the post trunc.
time column of \Cref{tab:rail_timing}.
Comparing the different the timings of different methods, the trade-off between
the block approach and the tangential variants becomes clear.
The increased computational costs in \tadi{} and \taditrunc{} compared to the
block methods comes from the additional linear systems of equations that have
to be solved in every step.
For this example, it is most efficient to simply perform a post truncation of the
classical block ADI method.

\begin{figure}[t]
  \centering
  \tikzexternalenable%
  \tikzsetnextfilename{rail_example}%
  \begin{tikzpicture}[font = \plotfontsize\normalfont]

  \pgfplotstableread{graphics/data/bilinear_step1_residual_adi.dat}\tableADI
  \pgfplotstableread{graphics/data/bilinear_step1_residual_eigtadi.dat}\tableEigTADI
\pgfplotstableread{graphics/data/bilinear_step1_residual_adi_trunc.dat}\tableADItrunc  \pgfplotstableread{graphics/data/bilinear_step1_residual_eigtadi_trunc.dat}\tableEigTADItrunc
  
  \begin{semilogyaxis}[%
    scale only axis,
    width                 = .8\linewidth,
    height                = .2\linewidth,
    xmin                  = 0,
    xmax                  = 400,
    ymin                  = 1e-12,
    ymax                  = 1e+1,
    ytick                 = {1e+0, 1e-4, 1e-8, 1e-12},
    xlabel                = {dimension of solution factor $\bL_{j}$},
    ylabel                = {normalized implicit residual},
    ylabel style          = {yshift = -.3em},
    scaled x ticks        = false,
    x tick label style    = {/pgf/number format/1000 sep={\,}},
    y tick label style    = {/pgf/number format/1000 sep={\,}},
    grid                  = major,
    legend cell align     = left,
    legend style          = {
      at        = {(1, 1)},
      anchor    = north east,
      align     = left}
  ]
  
    \addplot[adi_imp_style] 
      table[x = x_iter_rank, y = y_implicit] {\tableADI};
    \addlegendentry{\adi{}}

    \addplot[adi_trunc_imp_style, mark repeat = 8, mark phase = 2]
      table[x = x_iter_rank, y = y_implicit] {\tableADItrunc};
    \addlegendentry{\aditrunc{}}
    
    \addplot[eigtadi_imp_style, mark repeat = 25, mark phase = 0] 
      table[x = x_iter_rank, y = y_implicit] {\tableEigTADI};
    \addlegendentry{\tadi{}}

    \addplot[eigtadi_trunc_imp_style, mark repeat = 8, mark phase = 4]
      table[x = x_iter_rank, y = y_implicit] {\tableEigTADItrunc};
    \addlegendentry{\taditrunc{}}
    
  \end{semilogyaxis}
\end{tikzpicture}%
  \tikzexternaldisable%

  \caption{Convergence plot for the linear steel profile example:
    All methods converge fast to suitably low implicit ADI residuals.
    The tangential approaches \tadi{} and \taditrunc{} provide consistently
    smaller solution factors for the same level of accuracy compared to the block
    variants \adi{} and \aditrunc{}.}
  \label{fig:rail}
\end{figure}

\begin{table}[t]
  \centering
  \caption{Computation times for the linear steel profile example:
    The columns show the computations times of the respective solvers, the time
    for a post truncation step to bring all solution factors to the same size,
    and the total timings.
    Standard \adi{} is exceptionally fast compared to the other methods, while the
    internal truncation and tangential directions introduce a noticeable
    computational overhead.}
  \label{tab:rail_timing}

  \pgfkeys{/pgf/number format/1000 sep = {\,}}
  \begin{filecontents*}{graphics/data/bilinear_rail_step1_timing.dat}
Method	Solver_Time_s	Final_Trunc_s	Total_Time_s
ADI	2.815931	0.396410	3.212340
ADI (Int. Trunc)	9.371363	0.185982	9.557345
EigTADI (Int. Trunc)	16.674100	0.183618	16.857719
EigTADI	10.507982	0.357897	10.865879
\end{filecontents*}
\pgfplotstabletypeset[
    col sep               = tab,
    columns               = {Method, Solver_Time_s, Final_Trunc_s, Total_Time_s},
    every head row/.style = {before row = \toprule, after row = \midrule},
    every last row/.style = {after row = \bottomrule},
    font                  = \plotfontsize\normalfont,
    columns/Method/.style = {
      string type,
      column type    = l, 
      column name    = \textbf{method},
      string replace = {ADI}{\adi{}},
      string replace = {ADI (Int. Trunc)}{\aditrunc{}},
      string replace = {EigTADI}{\tadi{}},
      string replace = {EigTADI (Int. Trunc)}{\taditrunc{}} 
      },
    columns/{Solver_Time_s}/.style = {
      fixed,
      fixed zerofill,
      precision   = 4,
      column type = r,
      column name = {\textbf{solver time (s)}}
    },
    columns/{Final_Trunc_s}/.style = {
      fixed,
      fixed zerofill,
      precision   = 4,
      column type = r,
      column name = {\textbf{post trunc. time (s)}}
    },
    columns/{Total_Time_s}/.style = {
      fixed,
      fixed zerofill,
      precision   = 4,
      column type = r,
      column name = {\textbf{total time (s)}}
    },
    sort     = true,
    sort key = {Method},
    sort cmp = {string <}
  ]{graphics/data/bilinear_rail_step1_timing.dat}
\end{table}


\subsubsection{Simplified nonlinear heat transfer problem}%
\label{sec:bilinrail}

This problem is a modification of the previous example modeling the cooling
process of a steel bar.
The original model developed in~\cite{Saa03} naturally leads to a nonlinear
dynamical system of the form
\begin{equation} \label{eqn:bilinrail}
  \widetilde{\bE} \dot{\bx}(t) = \widetilde{\bA} \bx(t) +
    \sum\limits_{k = 1}^{m} \widetilde{\bN}_{k} \bx(t) \bu_{k}(t) +
    \widetilde{\bB} \bu(t), \quad
  \by(t) = \widetilde{\bC} \bx(t),
\end{equation}
where the bilinear terms are described by the matrices
$\widetilde{\bN}_{1}, \ldots, \widetilde{\bN}_{m} \in \R^{n \times n}$.
The data for this example can be found in~\cite{SaaKB23} and the matrix
dimensions are as in \Cref{sec:rail}.
In the design of low-dimensional models, systems analysis and controller design,
systems of the form~\cref{eqn:bilinrail}, multi-term Lyapunov equations
need to be solved.
A popular approach is the solution of such multi-term equations via the
solution of coupled Lyapunov equations of the form
\begin{equation} \label{eqn:bilyap}
  \bA \bX_{j} \bE^{\trans} + \bE \bX_{j} \bA^{\trans} +
    \sum\limits_{k = 1}^{m} \bN_{k} \bX_{j - 1} \bN_{k}^{\trans} +
    \bB \bR \bB^{\trans} = 0, 
\end{equation}
with the initial iterate $\bX_{0} = 0$;
see, for example,~\cite{ShaSS15}.
For the purpose as a numerical example in this work, we consider only
a single coupled Lyapunov equation of the form~\cref{eqn:bilyap},
namely for $j = 2$.
Furthermore, we note that while the iteration scheme described in~\cite{ShaSS15}
for the solution of multi-term Lyapunov equations applies an aggressive
rank truncation based on the current residual, this is not done here as we only
apply a single step of the approach and aim to illustrate the performance of
the tangential ADI in the setting of high-rank constant terms.
Given the low-rank approximation of the first non-zero Lyapunov solution
$\bX_{1} \approx \bL \bD \bL^{\trans}$, we can reformulate~\cref{eqn:bilyap}
into a classical indefinite Lyapunov equation with
\begin{equation} \label{eqn:modbilyap}
  \bA \bX_{2} \bE^{\trans} + \bE \bX_{2} \bA^{\trans} +
    \widetilde{\bB} \widetilde{\bR} \widetilde{\bB}^{\trans} = 0,
\end{equation}
where
\begin{equation*}
  \widetilde{\bB} = 
    \begin{bmatrix} \bB & \bN_{1} \bL & \ldots & \bN_{m} \bL \end{bmatrix}
  \quad\text{and}\quad
  \widetilde{\bR} =  \begin{bmatrix} \bR & & & & \\ & \bD & & & \\
    & & \ddots & \\ & & & & \bD \end{bmatrix}.
\end{equation*}
We note that the constant term in~\cref{eqn:modbilyap} has generically a large
rank due to the concatenation of previous solution factors, yet the solution
of~\cref{eqn:modbilyap} is typically numerically low rank.
For our numerical experiments, we use the rank-truncated version
of the solution to the Lyapunov equation (down to rank $209$) obtained for the
linear steel profile example in \Cref{sec:rail} and restrict to a single
nonzero $\bN_{k} = \widetilde{\bN}_{k}^{\trans}$ term.
The latter adjustment limits the size of the right-hand side to a level, which
allows us to apply the classical \adi{} method for the comparison.

The convergence results for this experiment are shown in \Cref{fig:bilinrail}.
The tangential approaches, \tadi{} and \taditrunc{}, and the block approach
with intermediate compression, \aditrunc{}, all effectively produce low-rank
solutions to~\cref{eqn:modbilyap} of rank about $300$.
This stands in stark contrast to the classical ADI, which accumulates about
$13\,500$ columns in outs solution factor approximation during execution.
The corresponding computation times are shown in \Cref{tab:bilinrail_timing}.
In this example, the tangential approaches strongly outperform the block ADI
variants with the fastest method being \tadi{} with a post truncation step to
further reduce the rank of the final solution approximation.
The \aditrunc{} approach uses significant amounts of computation time for the
intermediate truncations but allows to keep the size of the solution approximation
in a similar range compared to the tangential methods.
The classical block \adi{} runs faster than \aditrunc{} but underperforms in
terms of the total time due to the costs of the post truncation to reduce
the full solution approximation down to a reasonable size.

\begin{figure}[t]
  \centering
  \tikzexternalenable%
  \tikzsetnextfilename{bilinrail_example}%
  \begin{tikzpicture}[font = \plotfontsize\normalfont]
  \pgfplotstableread{graphics/data/bilinear_step2_residual_adi.dat}\tableADI
  \pgfplotstableread{graphics/data/bilinear_step2_residual_eigtadi.dat}\tableEigTADI
  \pgfplotstableread{graphics/data/bilinear_step2_residual_adi_trunc.dat}\tableADItrunc
  \pgfplotstableread{graphics/data/bilinear_step2_residual_eigtadi_trunc.dat}\tableEigTADItrunc
  
  \begin{loglogaxis}[%
    scale only axis,
    width                 = .8\linewidth,
    height                = .2\linewidth,
    xmin                  = 1,
    xmax                  = 16000,
    ymin                  = 1e-12,
    ymax                  = 1e+1,
    ytick                 = {1e+0, 1e-4, 1e-8, 1e-12},
    xlabel                = {dimension of solution factor $\bL_{j}$},
    ylabel                = {normalized implicit residual},
    ylabel style          = {yshift = -.3em},
    scaled x ticks        = false,
    x tick label style    = {/pgf/number format/1000 sep={\,}},
    y tick label style    = {/pgf/number format/1000 sep={\,}},
    grid                  = major,
    legend cell align     = left,
    legend style          = {
      at        = {(0, 0)},
      anchor    = south west,
      align     = left}
  ]
  
    \addplot[adi_imp_style] 
      table[x = x_iter_rank, y = y_implicit] {\tableADI};
    \addlegendentry{\adi{}}

    \addplot[adi_trunc_imp_style, mark repeat = 7, mark phase = 0]
      table[x = x_iter_rank, y = y_implicit] {\tableADItrunc};
    \addlegendentry{\aditrunc{}}
    
    \addplot[eigtadi_imp_style, mark repeat = 20, mark phase = 0]
      table[x = x_iter_rank, y = y_implicit] {\tableEigTADI};
    \addlegendentry{\tadi{}}
    
    \addplot[eigtadi_trunc_imp_style, mark repeat = 7, mark phase = 3]
      table[x = x_iter_rank, y = y_implicit] {\tableEigTADItrunc};
    \addlegendentry{\taditrunc{}}
    
  \end{loglogaxis}
\end{tikzpicture}%
  \tikzexternaldisable%

  \caption{Convergence plot for the simplified nonlinear rail example:
    The approximations computed by the tangential approaches \tadi{} and
    \taditrunc{} as well as truncated block method \aditrunc{} demonstrate
    the potential of efficiently compressing the solution approximations and saving
    memory costs during the iterations compared to the classical block approach.
    The roughly $13\,500$ columns in the solution accumulated in \adi{} are
    compressed by the other approaches during execution to about $300$.}
  \label{fig:bilinrail}
\end{figure}

\begin{table}[t]
  \centering
  \caption{Computation times for the simplified nonlinear rail example:
    The columns show the computations times of the respective solvers, the time
    for a post truncation step to bring all solution factors to the same size,
    and the total timings.
    Standard \adi{} outperforms \aditrunc{} in terms of solver runtime but
    require a significantly higher time to perform an additional post truncation
    step to reduce the size of the solution approximation.
    The tangential methods clearly outperform the block approaches, with \tadi{}
    being the fastest taking only $97.2$\,s including the post truncation step.}
  \label{tab:bilinrail_timing}

  \pgfkeys{/pgf/number format/1000 sep = {\,}}
  \begin{filecontents*}{graphics/data/bilinear_rail_step2_timing.dat}
Method	Solver_Time_s	Final_Trunc_s	Total_Time_s
ADI	138.578446	455.712107	594.290553
ADI (Int. Trunc)	201.304782	0.723513	202.028296
EigTADI (Int. Trunc)	110.718935	0.280343	110.999278
EigTADI	96.585020	0.652964	97.237983
\end{filecontents*}
\pgfplotstabletypeset[
    col sep               = tab,
    columns               = {Method, Solver_Time_s, Final_Trunc_s, Total_Time_s},
    every head row/.style = {before row = \toprule, after row = \midrule},
    every last row/.style = {after row = \bottomrule},
    font                  = \plotfontsize\normalfont,
    columns/Method/.style = {
      string type,
      column type    = l, 
      column name    = \textbf{method},
      string replace = {ADI}{\adi{}},
      string replace = {ADI (Int. Trunc)}{\aditrunc{}},
      string replace = {EigTADI}{\tadi{}},
      string replace = {EigTADI (Int. Trunc)}{\taditrunc{}} 
      },
    columns/{Solver_Time_s}/.style = {
      fixed,
      fixed zerofill,
      precision   = 4,
      column type = r,
      column name = {\textbf{solver time (s)}}
    },
    columns/{Final_Trunc_s}/.style = {
      fixed,
      fixed zerofill,
      precision   = 4,
      column type = r,
      column name = {\textbf{post trunc. time (s)}}
    },
    columns/{Total_Time_s}/.style = {
      fixed,
      fixed zerofill,
      precision   = 4,
      column type = r,
      column name = {\textbf{total time (s)}}
    },
    sort     = true,
    sort key = {Method},
    sort cmp = {string <}
  ]{graphics/data/bilinear_rail_step2_timing.dat}
\end{table}


\subsubsection{Vibrational response of a mechanical microstructure}%
\label{sec:butterfly}

The final numerical example comes from the modeling of the vibrational response
of a micro-mechanical gyroscope~\cite{morwiki_gyro, Bil05}.
The corresponding dynamical system is described by second-order differential
equations.
After reformulation into first-order form, we consider a Lyapunov equation
of the form~\cref{eqn:lyap_real}, where the coefficient matrices have the
following block structure
\begin{equation*}
  \bA = \begin{bmatrix} 0 & -\widetilde{\bK} \\ \eye{q} & -\widetilde{\bD}
    \end{bmatrix}, \quad
  \bE = \begin{bmatrix} \eye{q} & 0 \\ 0 & \widetilde{\bM} \end{bmatrix}, \quad
  \bB = \begin{bmatrix} \widetilde{\bC}^{\trans} \\ 0 \end{bmatrix},
\end{equation*}
where $\widetilde{\bM}, \widetilde{\bD}, \widetilde{\bK} \in \R^{q \times q}$
are symmetric positive definite.
For the center term $\bR$, we model again an indefinite matrix which represents
the performance of the system and the observation of perturbations.
The dimensions of the Lyapunov equation are $n =34\,722$ and $m = 12$.

While both methods compute suitable solution factors, the one constructed
by the tangential approach is about $20$\% smaller than the one
computed by the classical block approach.
In relation to the overall size of the problem ($n =34\,722$), such a reduction
yields significant improvements in terms of memory consumption.

The convergence and performance results for this example are shown in \Cref{fig:butterfly} and \Cref{tab:butterfly_timing}.
The implicit residual shown in \Cref{fig:butterfly} indicates significant
storage reductions for both the block and tangential method when intermediate
truncations are applied.
However, recomputation of the Lyapunov residuals at the end reveals
that \aditrunc{} and \taditrunc{} only achieve about
$9.3248 \times 10^{-6}$ and $1.2575 \times 10^{-5}$ for the true normalized
residuals, while for the other two approaches \adi{} and \tadi{} the implicit
and true residuals coincide.
These results can be related to the theory developed in~\cite{KueF20} implying
that early truncations may result in loss of accuracy.
The computations times reported in \Cref{tab:butterfly_timing} illustrate
again how the computational costs shift between \adi{} and \tadi{} from
memory requirements to the solution of additional linear systems.

\begin{figure}[t]
  \centering
  \tikzexternalenable%
  \tikzsetnextfilename{butterfly_example}%
  \begin{tikzpicture}[font = \plotfontsize\normalfont]

  \pgfplotstableread{graphics/data/butterfly_residual_adi.dat}\tableADI
  \pgfplotstableread{graphics/data/butterfly_residual_eigtadi.dat}\tableEigTADI
  \pgfplotstableread{graphics/data/butterfly_residual_adi_trunc.dat}\tableADItrunc
  \pgfplotstableread{graphics/data/butterfly_residual_eigtadi_trunc.dat}\tableEigTADItrunc
  
  \begin{semilogyaxis}[%
    scale only axis,
    width                 = .8\linewidth,
    height                = .2\linewidth,
    xmin                  = 0,
    xmax                  = 6500,
    ymin                  = 1e-12,
    ymax                  = 1e+1,
    xtick                 = {0, 1000, 2000, 3000, 4000, 5000, 6000},
    ytick                 = {1e+0, 1e-4, 1e-8, 1e-12},
    xlabel                = {dimension of solution factor $\bL_{j}$},
    ylabel                = {normalized implicit residual},
    ylabel style          = {yshift = -.3em},
    scaled x ticks        = false,
    x tick label style    = {/pgf/number format/1000 sep={\,}},
    y tick label style    = {/pgf/number format/1000 sep={\,}},
    grid                  = major,
    legend cell align     = left,
    legend style          = {
      at        = {(1, 1)},
      anchor    = north east,
      align     = left}
  ]
    
    \addplot[adi_imp_style] 
      table[x = x_iter_rank, y = y_implicit] {\tableADI};
    \addlegendentry{\adi{}}

    \addplot[adi_trunc_imp_style, mark repeat = 20, mark phase = 5] 
      table[x = x_iter_rank, y = y_implicit] {\tableADItrunc};
    \addlegendentry{\aditrunc{}}
    
    \addplot[eigtadi_imp_style, mark repeat = 250, mark phase = 0] 
      table[x = x_iter_rank, y = y_implicit] {\tableEigTADI};
    \addlegendentry{\tadi{}}
    
    \addplot[eigtadi_trunc_imp_style, mark repeat = 30, mark phase = 10] 
      table[x = x_iter_rank, y = y_implicit] {\tableEigTADItrunc};
    \addlegendentry{\taditrunc{}}
    
  \end{semilogyaxis}
\end{tikzpicture}%
  \tikzexternaldisable%

  \caption{Convergence plot for the micro-mechanical gyroscope
    example:
    The \tadi{} method computes a solution factor that is about 20\% smaller
    than in the classical \adi{}.
    The best memory usage and on-the-fly compression is done by
    \aditrunc{} and \taditrunc{} but it has to be noted that this results in
    this example in a significant loss of accuracy for the solution.}
  \label{fig:butterfly}
\end{figure}

\begin{table}[t]
  \centering
  \caption{Computation times for the micro-mechanical gyroscope
    example:
    The columns show the computations times of the respective solvers, the time
    for a post truncation step to bring all solution factors to the same size,
    and the total timings.
    By a slight margin, the fastest performing method is the classical \adi{} followed
    by the \aditrunc{}.
    The computation times for the tangential approaches are about five times
    larger illustrating the shift of computational costs away from memory
    towards the solution of linear systems of equations.}
  \label{tab:butterfly_timing}

  \pgfkeys{/pgf/number format/1000 sep = {\,}}
  \begin{filecontents*}{graphics/data/butterfly_gyro_timing.dat}
Method	Solver_Time_s	Final_Trunc_s	Total_Time_s
ADI	349.662495	119.327392	468.989887
ADI (Int. Trunc)	470.484153	0.494318	470.978471
EigTADI	2477.199156	87.313866	2564.513023
EigTADI (Int. Trunc)	2684.698978	0.891214	2685.590192
\end{filecontents*}
\pgfplotstabletypeset[
    col sep               = tab,
    columns               = {Method, Solver_Time_s, Final_Trunc_s, Total_Time_s},
    every head row/.style = {before row = \toprule, after row = \midrule},
    every last row/.style = {after row = \bottomrule},
    font                  = \plotfontsize\normalfont,
    columns/Method/.style = {
      string type,
      column type    = l, 
      column name    = \textbf{method},
      string replace = {ADI}{\adi{}},
      string replace = {ADI (Int. Trunc)}{\aditrunc{}},
      string replace = {EigTADI}{\tadi{}},
      string replace = {EigTADI (Int. Trunc)}{\taditrunc{}} 
      },
    columns/{Solver_Time_s}/.style = {
      fixed,
      fixed zerofill,
      precision   = 4,
      column type = r,
      column name = {\textbf{solver time (s)}}
    },
    columns/{Final_Trunc_s}/.style = {
      fixed,
      fixed zerofill,
      precision   = 4,
      column type = r,
      column name = {\textbf{post trunc. time (s)}}
    },
    columns/{Total_Time_s}/.style = {
      fixed,
      fixed zerofill,
      precision   = 4,
      column type = r,
      column name = {\textbf{total time (s)}}
    },
    sort     = true,
    sort key = {Method},
    sort cmp = {string <}
  ]{graphics/data/butterfly_gyro_timing.dat}
\end{table}


\section{Conclusions}%
\label{sec:conclusions}

In this work, we presented a tangential variant of the ADI iteration for the
construction of low-rank solutions to Lyapunov equations with indefinite
constant terms.
We derived the theoretical foundations of the method in complex arithmetic and
provided a variant of the theory as well as the algorithm that allows the
construction of real solution factors in the case that the coefficient matrices
of the Lyapunov equation are real.
An adaptive projection-based shift computation method and a heuristic
selection strategy for tangential directions have been presented to ensure
efficient computations with guaranteed convergence.
The numerical examples illustrate that the proposed tangential
approach together with the adaptive parameter selection
provide solution approximations of smaller size at the same
level of accuracy given by the classical ADI method.

While we motivated the restriction of the tangential directions to the
eigenvectors of the center matrix in the constant term of the Lyapunov
equation, it is not proven that there exists no other set of potential
tangential directions that also allows for the tangential ADI method to
converge efficiently.
In particular, further investigations are needed into why the mixing of
positive and negative eigenspaces for the selection of the tangential directions
leads to the divergence of the tangential ADI method.
Theoretical results on this topic will have strong implications on the
application of inexact and compressed Krylov subspaces in other research areas
such as model order reduction or rational function interpolation.


\section*{Declaration of Competing Interest}

The authors declare that they have no known competing financial interests or
personal relationships that could have appeared to influence the work reported
in this paper.


\section*{Data Availability}

The source codes, data, and results of the numerical experiments are available
at~\cite{supSmiW26}.


\addcontentsline{toc}{section}{References}
\bibliographystyle{plainurl}
\bibliography{bibtex/myref}

\end{document}